\documentclass[number,citesort,dvips]{arxbj}

\aid{0}
\volume{17}
\issue{4}
\pubyear{2011}
\firstpage{1217}
\lastpage{1247}
\doi{10.3150/10-BEJ312}

\makeatletter
\newtheorem{theorem}{Theorem}%[section]

\newtheorem{proposition}{Proposition}
\newtheorem{corollary}{Corollary}
\newtheorem{lemma}{Lemma}
\newproclaim{nnot}{Notation}

\newcommand{\E}{\mathbb{E}}
\newcommand{\N}{\mathbb{N}}
\newcommand{\R}{\mathbb{R}}
\newcommand{\Z}{\mathbb{Z}}

\renewcommand{\P}{\mathbb{P}}
\newcommand{\be}{\mathbf{e}}
\newcommand{\FF}{\mathcal{F}}
\renewcommand{\d}{\mathrm{d}}
\newcommand{\eps}{\varepsilon}
\newcommand{\ind}{\mathbf{1}}
\makeatother

\begin{document}
\begin{frontmatter}

\title{Self-similar scaling limits of non-increasing Markov chains}
\runtitle{Self-similar scaling limits of non-increasing Markov chains}

\begin{aug}
\author{\fnms{B\'{e}n\'{e}dicte} \snm{Haas}\ead[label=e1]{haas@ceremade.dauphine.fr}\corref{}}
\and
\author{\fnms{Gr\'{e}gory} \snm{Miermont}\ead[label=e2]{gregory.miermont@math.u-psud.fr}}

\runauthor{B. Haas and G. Miermont}
\address{Universit\'{e} Paris-Dauphine, Paris, France. \printead{e1}}
\address{Universit\'{e} Paris-Sud, Paris, France. \printead{e2}}
\end{aug}

% HISTORY:
\received{\smonth{9} \syear{2009}}
\revised{\smonth{6} \syear{2010}}

% ABSTRACT
%
\begin{abstract}
We study scaling limits of non-increasing Markov chains with values in
the set of non-negative integers, under the assumption that the large jump
events are rare and happen at rates that behave like a negative power
of the current state. We show that the chain starting from $n$ and
appropriately rescaled, converges in distribution, as $n \rightarrow
\infty$, to a non-increasing self-similar Markov process. This
convergence holds jointly with that of the rescaled absorption time to
the time at which the self-similar Markov process reaches first~0.

We discuss various applications to the study of random walks with a
barrier, of the number of collisions in $\Lambda$-coalescents that do
not descend from infinity and of non-consistent regenerative
compositions. Further applications to the scaling limits of Markov
branching trees are developed in our paper, Scaling limits of {M}arkov
branching trees,
with applications to Galton--Watson and random unordered trees (2010).
\end{abstract}

% KEYWORDS
%
\begin{keyword}
\kwd{absorption time}
\kwd{$\Lambda$-coalescents}
\kwd{random walks with a barrier}
\kwd{regenerative compositions}
\kwd{regular variation}
\kwd{self-similar Markov processes}
\end{keyword}

\end{frontmatter}

%s1 ###
\section{Introduction and main results}
\label{Intro}
%%%%%%%%%%%%%%%%%%%

Consider a Markov chain taking values in the set of non-negative
integers $\Z_+=\{0,1,2,\ldots\}$, and with non-increasing paths. We are
interested in the asymptotic behavior in distribution of the chain
started from $n$, as $n$ tends to $\infty$. Our main assumption is
(roughly speaking) that the chain, when in state $n$, has a ``small''
probability, of order $c_\eps n^{-\gamma}$ for some $\gamma>0$ and some
$c_\eps>0$, of accomplishing a negative jump with size in
$[n\varepsilon,n]$, where $0<\eps<1$. A typical example is constructed
from a random walk $(S_k,k \geq0)$ with non-negative steps with tail
distribution proportional to $n^{-\gamma}$ as $n$ tends to $\infty$,
for some $\gamma\in(0,1)$, by considering the Markov chain starting
from $n$: $(\max(n-S_k,0), k\geq0)$. An explicit example is provided by
the step distribution $q_n=(-1)^{n-1}{{\gamma}\choose{n}}$, $n \geq1$.

Under this main assumption, we show in Theorem \ref{CvChain} that the
chain started from $n$, and properly rescaled in space and time,
converges in distribution in the Skorokhod space to a non-increasing
self-similar Markov process. These processes were introduced and
studied by Lamperti \cite{lamperti62,Lamperticaract}, under the name of
\textit{semi-stable} processes, and by many authors since then. Note that
Stone \cite{stone63} discusses limit theorems for birth-and-death
chains and diffusions that involve self-similar Markov processes, but
in a context that is very different from ours.

A quantity of particular interest is the absorption time of the chain,
that is, the first time after which the chain remains constant. We show
in Theorem \ref{cvdeathtimes} that jointly with the convergence of
Theorem \ref{CvChain}, the properly rescaled absorption time converges
to the first time the limiting self-similar Markov process hits $0$. In
fact, we even show that all positive moments of the rescaled absorption
time converge.

These results have applications to a number of problems considered in
the literature, such as the random walk with a barrier \cite{IM08} when
the step distribution is in the domain of attraction of a stable random
variable with index in $(0,1)$, or the number of coalescing events in a
$\Lambda$-coalescent that does not come down from infinity
\cite{GIM08,IM08}. It also allows us to recover some results by Gnedin,
Pitman and Yor \cite{GPYalpha} for the number of blocks in regenerative
composition structures, and to extend this result to the case of
``non-consistent compositions''. One of the main motivations for the
present study was to provide a unified framework to treat such
problems, which can all be translated in terms of absorption times of
non-increasing Markov chains. Moreover, the convergence of the rescaled
Markov chain as a process, besides the convergence of the absorption
time, provides new insights on these results. Finally, our main results
are also a starting point for obtaining the scaling limits of a large
class of random trees satisfying a simple discrete branching property.
This is the object of the paper \cite{HMMarkovBranching}.

Let us now present our main results and applications in a more formal
way. Implicitly, all the random variables in this paper are defined on
a probability space $(\Omega,\FF,\P)$.

\begin{nnot*}
For two positive sequences $x_n,y_n,n\geq0$,
the notation $x_n\sim y_n$ means that $x_n/y_n$ converges to $1$ as
$n\to\infty$.
\end{nnot*}

%%%%%%%%%%%%%%%%%%%%%%%%%%%%%%%
%s1.1 ###
\subsection{Scaling limits of non-increasing Markov chains}
\label{secscal}
%%%%%%%%%%%%%%%%%%%%%%%%%%%%%%%

For every $n\geq0$, consider a non-negative sequence $(p_{n,k},0\leq
k\leq n)$ that sums to $1$,
\[
\sum_{k=0}^np_{n,k}=1 .
\]
We view the latter as a probability
distribution on $\{0,1,\ldots,n\}$, and view the family $(p_{n,k},\break0
\leq k\leq n)$ as the transition probabilities for a discrete-time
Markov chain, which takes integer values and has non-increasing paths.
We will denote by $(X_n(k),k\geq0)$ such a Markov chain, starting at
the state $X_n(0) = n$. For every $n\geq1$, we let
%define a probability distribution
$p^*_n$ be the law on $[0,1]$ of $X_n(1)/n$, so that
%as the push-forward of $(p_{n,k},0\leq k\leq n)$ under
%$x\mapsto x/n$:
%
\[
p^*_n(\d{x})=\sum_{k=0}^{n}p_{n,k}\delta_{k/n}(\d{x}) .
\]
Our main
assumption all throughout the paper will be the following hypothesis.

\begin{h*}
There exist:
\begin{itemize}
\item a sequence $(a_n,n\geq1)$ of the form $a_n =
n^{\gamma}\ell(n)$, where $\gamma>0$ and $\ell\dvtx \R_+\to(0,\infty
)$ is a
function that is slowly varying at $\infty$,
\item a non-zero, finite,
non-negative measure $\mu$ on $[0,1]$,
\end{itemize}
such that the following weak convergence of finite measures on $[0,1]$
holds:
%
%e1 ###
\begin{equation}
\label{Hyp} a_n(1-x) p^*_n(\d{x})
\mathop{\stackrel{(w)}{\longrightarrow}}_{n\to\infty}\mu(\d{x}) .
\end{equation}
\end{h*}

This means that a jump of the process $X_n/n$ from $1$ to $x\in(0,1)$
occurs with a small intensity $a_n^{-1}\mu(\d{x})/(1-x)$, and indicates
that an interesting scaling limit for the Markov chain $X_n$ should
arise when rescaling space by $n$ and time by $a_n$. Also, note that
$\mu([0,1])/a_n$ is equivalent as $n\to\infty$ to the expectation of
the first jump of the chain $X_n/n$, and this converges to 0 as $n
\rightarrow\infty$. The role of the factor $(1-x)$ in (\ref{Hyp}) is
to temper the contribution of very small jumps in order to evaluate the
contribution of larger jumps.

Of course, in (H), the sequence $a=(a_n,n\geq1)$, the function
$\ell$ and the measure $\mu$ are not uniquely determined. One can
simultaneously replace $a$ by $ca$ and $\mu$ by $c\mu$ for any given
$c>0$. Also, one can replace $\ell$ by any function that is equivalent
to it at infinity. However, it is clear that $\mu$ is determined up to
a positive multiplicative constant (with a simultaneous change of the
sequence $a$ as depicted above), and that $\gamma$ is uniquely
determined.

We will soon see that hypothesis (H) appears very naturally in
various situations. It is also very general, in the sense that there
are no restrictions on the sequences $(a_n,n\geq1)$ or measures $\mu$
that can arise. Here is a formal statement, which is proved at the end
of Section \ref{Markovchain}.

\begin{proposition}\label{sec:scaling-limits-non-3}
For any finite measure $\mu$ on $[0,1]$ and any sequence of the form
$a_n = n^{\gamma}\ell(n)$ where $\gamma>0$ and $\ell\dvtx \R_+\to
(0,\infty)$ is slowly varying at $\infty$, one can find a sequence of
probability vectors $((p_{n,k}, 0 \leq k \leq n), n \geq0 )$ such that
\textup{(\ref{Hyp})} holds.
\end{proposition}

We now describe the objects that will arise as scaling limits of $X_n$.
For $\lambda>0 $ and $x\in[0,1)$, let
%
%e2 ###
\begin{equation}\label{eq:3}
[\lambda]_x=\frac{1-x^\lambda}{1-x} ,\qquad 0\leq x<1 ,
\end{equation}
and set $[\lambda]_1=\lambda$. For each $\lambda> 0$, this defines a
continuous function $x\mapsto[\lambda]_x$ on $[0,1]$. If $\mu$ is a
finite measure on $[0,1]$, then the function $\psi$ defined for
$\lambda>0$ by
%
%e3 ###
\begin{equation}\label{eq:1}
\psi(\lambda):=\int_{[0,1]}[\lambda]_x\mu(\d{x})
\end{equation}
and extended at $0$ by $\psi(0):=\lim_{\lambda\downarrow0}
\psi(\lambda)=\mu(\{ 0\})$ is the Laplace exponent of a subordinator.
To see this, let $\mathtt{k}=\mu(\{0\}), \mathtt{d}=\mu(\{1\})$, so
that $\psi$ can be written in the usual L{\'{e}}vy--Khintchine form:
\[
\psi(\lambda)=
\mathtt{k}+\mathtt{d}\lambda+\int_{(0,1)}(1-x^\lambda)\frac{\mu
(\d{x})}{1-x} =\mathtt{k}+\mathtt{d}\lambda+
\int_0^\infty(1-\mathrm{e}^{-\lambda y}) \omega(\mathrm d y) ,
\]
where $\omega$ is the push-forward of the measure
$(1-x)^{-1}\mu(\d{x})\ind_{\{0<x<1\}}$ by the mapping $x\mapsto
-\log
x$. Note that $\omega$ is a $\sigma$-finite measure on $(0,\infty)$
that integrates $y\mapsto y\wedge1$, as it ought. Conversely, any
Laplace exponent of a (possibly killed) subordinator can be put in the
form (\ref{eq:1}) for some finite measure $\mu$.

Now, let $\xi$ be a subordinator with Laplace exponent $\psi$. This
means that the process $(\xi_t,t\geq0)$ is a non-decreasing L\'{e}vy
process with
\[
\E[\exp(-\lambda\xi_t)]=\exp(-t\psi(\lambda)) ,\qquad t,\lambda\geq0 .
\]
Note in particular that the subordinator is
killed at rate $\mathtt{k}\geq0$. The function $t \in[0,\infty)
\rightarrow
\int_0^t \exp(-\gamma\xi_r)\,\mathrm d r$ is continuous,
non-decreasing and its limit at infinity, denoted by
\[
I:=\int_0^\infty\exp(-\gamma\xi_r)\,\d{r},
\]
is a.s. finite. Standard properties of this random variable are studied
in \cite{CarmonaPetitYor}. We let $\tau\dvtx [0,I)\to\R_+$ be its inverse
function, and set $\tau(t)=\infty$ for $t\geq I$. The process
%
%e4 ###
\begin{equation}\label{eq:2}
Y(t):=\exp\bigl(-\xi_{\tau(t)}\bigr) ,\qquad t\geq0,
\end{equation}
is a non-increasing self-similar Markov process starting from 1.
%Self-similar Markov processes have been studied by Lamperti
%process),
Recall from \cite{Lamperticaract} that if $\P_x$ is the law of an
$\R_+$-valued Markov process $(M_t,t\geq0)$ started from $M_0=x\geq
0$, then the process is called self-similar with exponent $\alpha>0$ if
the law of $(r^{-\alpha}M_{rt},t\geq0)$ under $\P_x$ is
$\P_{r^{-\alpha}x}$, for every $r>0$ and $x\geq0$.

In this paper, all processes that we consider belong to the space
$\mathcal{D}$ of c\`{a}dl\`{a}g, non-negative functions from
$[0,\infty)$
to $\R$. This space is endowed with the Skorokhod metric, which
makes it a Polish space. We refer to \cite{EK}, Chapter 3.5, for
background on the topic. We recall that $\lfloor r \rfloor$ denotes
the integer part of the real number $r$.

\begin{theorem}\label{CvChain}
For all $t \geq0$ and all $n \in\mathbb N$, we let
\[
Y_n(t):=\frac{X_n( \lfloor a_nt \rfloor)}{n}.
\]
Then, under the assumption \textup{(H)}, we have the following convergence
in distribution
\[
Y_n \mathop{\stackrel{(d)}{\longrightarrow}}_{n\to\infty} Y
\]
for the Skorokhod topology on $\mathcal{D}$, where $Y$ is defined at
(\ref{eq:2}).
\end{theorem}

A theorem by Lamperti \cite{Lamperticaract} shows that any c\`{a}dl\`{a}g,
non-increasing, non-negative, self-similar Markov process (started from
1) can be written in the form (\ref{eq:2}) for some subordinator $\xi$
and some $\gamma>0$. In view of this, Theorem \ref{CvChain}, combined
with Proposition~\ref{sec:scaling-limits-non-3}, implies that every
non-increasing, c\`{a}dl\`{a}g, self-similar Markov process is the weak scaling
limit of a non-decreasing Markov chain with rare large jumps.

In fact, as the proof of Theorem \ref{CvChain} will show, a more
precise result holds. With the above notations, for every $t\geq0$,
we let $Z(t)=\exp(-\xi_t)$, so that $Y(t)=Z(\tau(t))$. Let also
\[
\tau_n^{-1}(t)=\inf\biggl\{u\geq0\dvtx \int_0^{u} Y_n^{-\gamma}(r)
\,\mathrm d
r>t\biggr\} ,\qquad t\geq0 ,
\]
and $Z_n(t):=Y_n(\tau_n^{-1}(t))$.\vadjust{\goodbreak}

\begin{proposition}\label{sec:scaling-limits-non-2}
Under the same hypotheses and notations as Theorem \ref{CvChain}, one
has the joint convergence in distribution
\[
(Y_n,Z_n) \mathop{\stackrel{(d)}{\longrightarrow}}_{n\to\infty}
(Y, Z)
\]
for the product topology on $\mathcal{D}^2$.
\end{proposition}

%%%%%%%%%%%%%%%%%%%%%%%%%
%s1.2 ###
\subsection{Absorption times}
%%%%%%%%%%%%%%%%%%%%%%%%%

Let $\mathcal{A}$ be the set of absorbing states of the chain, that is,
\[
\mathcal{A}:=\{k\in\Z_+ \dvtx p_{k,k}=1 \} .
\]
Under assumption
(H) it is clear that $\mathcal{A}$ is finite, and not empty since
it contains at least $0$. It is also clear that the absorbing time
\[
A_n:=\inf\{k \in\mathbb Z_+ \dvtx X_{n}(k) \in\mathcal A\}
\]
is a.s. finite. For $(Y,Z)=(\exp(-\xi_{\tau}),\exp(-\xi))$ defined as
in the previous subsection, we let $\sigma=\inf\{t\geq0\dvtx Y(t)=0\}$. Then
it holds that
%
%e5 ###
\begin{equation}\label{eq:9}
\sigma=\int_0^\infty\exp(-\gamma\xi_r)\,\d{r} ,
\end{equation}
which is a general fact that we recall (\ref{sigma}) in Section
\ref{tc} below.

\begin{theorem}
\label{cvdeathtimes}
Assume \textup{(H)}. Then, as $n \rightarrow\infty$,
\[
\frac{A_n}{a_n}\stackrel{(d)}{\rightarrow}\sigma,
\]
and this holds jointly with the
convergence in law of $(Y_n,Z_n)$ to $(Y,Z)$ as stated in Proposition~\ref{sec:scaling-limits-non-2}. Moreover, for all $p \geq0$,
\[
\mathbb E\biggl[ \biggl(\frac{A_n}{a_n} \biggr)^p\biggr] \rightarrow\mathbb E[ \sigma^p] .
\]
When $p \in\mathbb{Z}_+$, the limiting moment $\E[\sigma^p]$ is
equal to $p ! / \prod_{i=1}^p \psi(\gamma i) $.
\end{theorem}

Note that even the first part of this result is not a direct
consequence of Theorem \ref{CvChain} since convergence of functions in
$\mathcal{D}$ does not lead, in general, to the convergence of their
absorption times (when they exist).

%%%%%%%%%%%%%%%%%%%%%%%%%%%%%%%%%%%
%s1.3 ###
\subsection{Organization of the paper}\label{sec:organisation-paper}
%%%%%%%%%%%%%%%%%%%%%%%%%%%%%%%%%%%

We start in Section \ref{sec:applications} with a series of
applications of Theorems \ref{CvChain} and \ref{cvdeathtimes} to random
walks with a barrier, $\Lambda$-coalescents and non-consistent
regenerative compositions. Most of the proofs of these results, as well
as further developments, are postponed to Sections \ref{RW} (for the
random walks with a barrier) and \ref{Coalescents} (for $\Lambda$-coalescents).

In the preliminary Section \ref{prem}, we gather some basic facts
needed for the proofs of Theorems~\ref{CvChain} and \ref{cvdeathtimes}
and Proposition~\ref{sec:scaling-limits-non-2}, which are undertaken in
Section \ref{Markovchain}. The proof of Proposition~\ref{sec:scaling-limits-non-2} and Theorem \ref{CvChain} will be
obtained by a classical two-step approach: first, we show that the laws
of $(Y_n,Z_n),n\geq1$ form a tight family of probability distributions
on $\mathcal{D}^2$. Then, we will show that the only possible limiting
distribution is that of $(Y,Z)$. This identification of the limit will
be obtained via a simple martingale problem. Tightness is studied in
Section \ref{Tightness} and the characterization of the limits in
Section \ref{Limit}. In both cases, we will work with some sequences of
martingales related to the chains $X_n$, which are introduced in
Section \ref{martingales}. The convergence of $(Y_n,Z_n)$ to $(Y,Z)$ is
a priori not sufficient to get the convergence of the absorption times,
as stated in Theorem \ref{cvdeathtimes}. This will be obtained in
Section \ref{Absorption}, by first showing that $t^{\beta}\mathbb
E[Z_n(t)^{\lambda} ]$ is uniformly bounded for every $\beta>0$.

Last, a proof of Proposition~\ref{sec:scaling-limits-non-3} is given at
the end of Section \ref{Markovchain}.

%%%%%%%%%%%%%%%%%%%%%%%%%%%
%s2 ###
\section{Applications}\label{sec:applications}
%%%%%%%%%%%%%%%%%%%%%%%%%%%

%s2.1 ###
\subsection{Random walk with a barrier}\label{sec:random-walk-with}

Let $q=(q_k,k\geq0)$ be a non-negative sequence with total sum
$\sum_kq_k=1$, which is interpreted as a probability distribution on
$\Z_+$. We assume that $q_0<1$ in order to avoid trivialities. For
$n\geq0$, we let
\[
\overline{q}_n=\sum_{k> n}q_k ,\qquad n\geq0 .
\]
%
%Let $(\zeta_i,i\geq
%0)$ be an i.i.d. sequence of random variables with distibution
%$q$.

The random walk with a barrier is a variant of the usual random walk
%$S_k=\sum_{i=1}^k\zeta_i, k\geq0$
%with increments
%$\zeta_1,\zeta_2,\ldots$.
with step distribution $q$. Informally, every step of the walk is
distributed as $q$, but conditioned on the event that it does not bring
the walk to a level higher than a given value $n$. More formally, for
every $n$, we define the random walk with barrier $n$ as the Markov
chain $(S_k^{(n)},k\geq0)$ starting at 0, with values in $\{0,1,2,\ldots,n\}$ and
%
% the random walk $S$ conditioned to remain
%below $n$ forever:
%$$S^{(n)}\overset{(d)}{=}S \mbox{ given } \{S_k\leq
%n\mbox{ for all }k\}\, .$$
%Of course, the conditioning is singular, so
%one has to be a little more careful, for instance by conditioning on
%the event that $\{S_k\leq n,0\leq k\leq M\}$ and letting
%$M\to\infty$. It is immediate to see that this does define a Markov
%process with values in $\{0,1,\ldots,n\}$ and
with transition probabilities
\[
q^{(n)}_{i,j}=\cases{
\displaystyle\frac{q_{j-i}}{1-\overline{q}_{n-i}},&\quad if $\overline{q}_{n-i}<1$,\cr
\ind_{\{j = i\}}, &\quad if $\overline{q}_{n-i}=1$,
}\qquad 0\leq i\leq j\leq n .
\]
(This definition is
not exactly the same as in
\cite{IM08}, but the absorption time $A_n$ is exactly the random
variable $M_n$, which is the main object of study in this paper. We
will comment further on this point in Section \ref{RW}.)

To explain the definition, note that when $\overline{q}_r<1$,
$(q_k/(1-\overline{q}_r),0\leq k\leq r)$ is the law of a random
variable with
distribution $q$, conditioned to be in $\{0,\ldots,r\}$. When
$\overline{q}_r=1$, the quotient is not well defined, and we choose the
convention that the conditioned law is the Dirac measure at $\{0\}$. In
other words, when the process arrives at a state $i$ such that
$\overline{q}_{n-i}=1$, so that every jump with distribution $q$ would be
larger than $n-i$, we choose to let the chain remain forever at state~$i$. Of course, the above discussion is not needed when $q_0>0$.

As a consequence of the definition, the process
\[
X_n(k) = n-S^{(n)}_k ,\qquad k\geq0,
\]
is a Markov process with
non-increasing paths, starting at $n$, and with transition
probabilities
%
%e6 ###
\begin{equation}
\label{eq:5} p_{i,j}=\frac{q_{i-j}}{1-\overline{q}_i} ,\qquad0\leq j\leq i ,
\end{equation}
with the convention that $p_{i,j}=\ind_{\{j=i\}}$ when $\overline{q}_i=1$.
The probabilities (\ref{eq:5}) do not depend on $n$, so this falls
under our basic framework. As before, we let $A_n$ be the absorbing
time for $X_n$.

\begin{theorem}
\label{ThRW}
\textup{(i)} Let $\gamma\in(0,1)$, and assume that
$\overline{q}_n
= n^{-\gamma}\ell(n),$ where $\gamma\in(0,1)$ and $\ell$ is slowly
varying at $\infty$. Let $\xi$ be a subordinator with Laplace exponent
\[
\psi(\lambda)= \int_0^\infty(1-\mathrm{e}^{-\lambda y}) \frac{\gamma
\mathrm{e}^{-y}\,\d{y}}{(1-\mathrm{e}^{-y})^{\gamma+1}} ,\qquad\lambda\geq0 ,
\]
and
let
\[
\tau(t)=\inf\biggl\{u\geq0\dvtx \int_0^u \exp(-\gamma\xi_r)\,\d{r}>t
\biggr\} ,\qquad
t\geq0 .
\]
Then,
\[
\biggl(\frac{X_n(\lfloor t/\overline{q}_n\rfloor)}{n}\biggr)
\mathop{\stackrel{(d)}{\longrightarrow}}_{n\to\infty}
\bigl(\exp\bigl(-\xi_{\tau(t)}\bigr),t\geq0\bigr) ,
\]
jointly with the convergence
\[
\overline{q}_nA_n\mathop{\stackrel{(d)}{\longrightarrow}}_{n\to\infty}
\int_0^\infty\exp(-\gamma\xi_t)\,\d{t} .
\]
For the latter, the
convergence of all positive moments also holds.

\textup{(ii)} Assume that $m:=\sum_{k=0}^{\infty}kq_k$ is finite. Then
\[
\biggl(\biggl(\frac{X_n(\lfloor tn\rfloor)}{n},t \geq
0\biggr),\frac{A_n}{n}\biggr)\mathop{\stackrel{(P)}{\longrightarrow}}_{n\to\infty}\bigl(\bigl( \bigl( (1-mt)\vee0\bigr), t \geq0 \bigr),1/ m\bigr) ,
\]
in probability in
$\mathcal{D}\times\R_+$. Convergence of all positive moments also
holds for the second components.
\end{theorem}

Of course, this will be proved by checking that (H) holds for
transition probabilities of the particular form~(\ref{eq:5}), under the
assumption of~\ref{ThRW}. This result encompasses Theorems 1.1 and 1.4
in \cite{IM08}. Note that Theorems 1.2 and 1.5 in the latter reference
give information about the deviation for $A_n$ around $n/m$ in case
(ii) of Theorem \ref{ThRW} above, under some assumptions on $q$ (saying
essentially that a random variable with law $q$ is in the domain of
attraction of a stable law with index in $[1,2]$, as opposed to $(0,1)$
in Theorem \ref{ThRW}). See also \cite{DeDhSJ08} for related results in
a different context.

The L\'{e}vy measure of the subordinator $\gamma\xi$ involved in
Theorem \ref{ThRW} is clearly given by
\[
\exp(-x/\gamma)\bigl(1-\exp(-x/\gamma)\bigr)^{-\gamma-1}\,\mathrm{d}x \,\mathbf
1_{\{ x
\geq0\}}.
\]
Bertoin and Yor \cite{BYFacExp} show that the variable $\int
_0^{\infty}
\exp(-\gamma\xi_r) \,\mathrm{d}r$ is then distributed as\break
$\Gamma(1-\gamma)^{-1}\tau_{\gamma}^{-\gamma}$, where $\tau
_{\gamma}$
is a stable random variable with Laplace transform $\mathbb E[
\exp(-\lambda\tau_{\gamma})]=\exp(-\lambda^{\gamma})$.

%s2.2 ###
\subsection{On collisions in $\Lambda$-coalescents that do
not come down from infinity}\label{sec:numb-coll-lambda}

We first briefly recall the definition and basic properties of a
$\Lambda$-coalescent, referring the interested reader to
\cite{jp99cmc,Sagitov} for more details.

Let $\Lambda$ be a finite measure on $[0,1]$. For $r\in\N$, a
$(\Lambda,r)$-coalescent is a Markov process $(\Pi_r(t),t\geq0)$
taking values in the set of partitions of $\{1,2,\ldots,r\}$, which is
monotone in the sense that $\Pi_r(t')$ is coarser than $\Pi_r(t)$ for
every $t'>t$. More precisely, $\Pi_r$ only evolves by steps that
consist of merging a certain number (at least $2$) of blocks of the
partition into one, the other blocks being left unchanged. Assuming
that $\Pi_r(0)$ has $n$ blocks, the rate of a collision event involving
$n-k+1$ blocks, bringing the process to a state with $k$ blocks, for
some $1 \leq k \leq n-1$, is given by
\[
g_{n,k}=\pmatrix{n\cr k-1} \int_{[0,1]} x^{n-k-1}(1-x)^{k-1}
\Lambda(\mathrm{d} x) ,
\]
and the blocks that intervene in the merging event are uniformly
selected among the ${{n}\choose{k-1}}$ possible choices of $n-k+1$ blocks
out of $n$. Note that these transition rates depend only on the number
of blocks present at the current stage. In particular, they do not
depend on the particular value of $r$.

A $\Lambda$-coalescent is a Markov process $(\Pi(t),t\geq0)$ with
values in the set of partitions of $\N$, such that for every $r\geq1$,
the restriction $(\Pi|_{[r]}(t),t\geq0)$ of the process to
$\{1,2,\ldots,r\}$ is a $(\Lambda,r)$-coalescent. The existence (and
uniqueness in law) of such a process is discussed in \cite{jp99cmc}.
The most celebrated example is the Kingman coalescent obtained for
$\Lambda=\delta_0$.

The $\Lambda$-coalescent $(\Pi(t),t\geq0)$ is said to \textit{come down
from infinity} if, given that $\Pi(0)=\{\{i\},i\geq1\}$ is the
partition of $\N$ that contains only singletons, $\Pi(t)$ a.s. has a
finite number of blocks for every $t>0$. When the coalescent does not
come down from infinity, it turns out that $\Pi(t)$ has a.s. infinitely
many blocks for every $t\geq0$, and we say that the coalescent \textit
{stays infinite}. See \cite{schw00infty} for more details and a nice
criterion for the property of coming down from infinity. By Lemma 25 in
\cite{jp99cmc}, the $\Lambda$-coalescent stays infinite if
$\int_{[0,1]}x^{-1}\Lambda(\d{x})<\infty$.

Starting with $n$ blocks in a $(\Lambda,r)$-coalescent (or in a
$\Lambda$-coalescent), let $X_n(k)$ be the number of blocks after $k$
coalescing events have taken place. Due to the above description, the
process $(X_n(k),k\geq0)$ is a Markov chain with transition
probabilities given by
%
%e7 ###
\begin{eqnarray}
\label{defprobacoal}
p_{n,k}&=&\mathbb P\bigl(X_n(1)=k
\bigr)=\frac{g_{n,k}}{g_n}\nonumber
\\
&=&\frac{1}{g_n}\pmatrix{n\cr k-1} \int_{[0,1]}
x^{n-k-1}(1-x)^{k-1} \Lambda(\mathrm{d} x),\qquad 1 \leq k \leq n-1,
\end{eqnarray}
where $g_n$ is the total transition rate $g_n =\sum_{k=1}^{n-1}
g_{n,k}$. This chain always gets absorbed at~$1$.

The total number of collisions in the coalescent coincides with the
absorption time $A_n := \inf\{k \dvtx X_n(k)=1 \}$. There have been many
studies on the asymptotic behavior of $A_n$ as $n \to\infty$
\cite{GIM08,GY07,IM08,IMM09}, in contexts that mostly differ from ours
(see the comments below). For $u\in(0,1]$ we let
%
%e8 ###
\begin{equation}\label{eq:6}
h(u)=\int_{[u,1]} x^{-2}\Lambda(\d{x}) .
\end{equation}
We are interested in cases where $\lim_{u\downarrow0}h(u)=\infty$ but
$\int_0^1x^{-1}\Lambda(\d{x})<\infty$, so the coalescent stays infinite
by the above discussion.

%knowledge, on the whole process of collisions $X_n$.
%
%. We discuss below the connexions with
%our result, which is a consequence of Theorems \ref{CvChain} and

\begin{theorem}
\label{theolambda}
Let $\gamma\in(0,1)$. We assume that the function
$h$ is regularly varying at 0 with index $-\gamma$. Let $\xi$ be a
subordinator with Laplace exponent
%
%e9 ###
\begin{equation}
\label{eq:7}
\psi(\lambda)=\frac{1}{\Gamma(2-\gamma)}\int_0^1\bigl(1-(1-x)^\lambda\bigr)
x^{-2}\Lambda(\d{x}) ,\qquad\lambda\geq0 ,
\end{equation}
and let
%
%e10 ###
\begin{equation}
\label{rho} \tau(t):=\inf\biggl\{u \geq0\dvtx \int_0^u \exp(-\gamma
\xi_r)
\,\d{r}
>t \biggr\} ,\qquad t\geq0 .
\end{equation}
Then,
%
%e11 ###
\begin{equation}
\label{coalambda}
\biggl( \frac{X_n(\lfloor h(1/n) t\rfloor)}{n}, t
\geq
0\biggr)\mathop{\stackrel{(d)}{\longrightarrow}}_{n\to\infty}
\exp(-\xi_{\tau}).
\end{equation}
Moreover, jointly with (\ref{coalambda}), it holds that
%
%e12 ###
\begin{equation}
\label{nblambda}
\frac{A_n}{h(1/n)}
\mathop{\stackrel{(d)}{\longrightarrow}}_{n\to\infty} \int
_0^{\infty}
\exp(-\gamma\xi_r) \,\mathrm{d}r,
\end{equation}
and there is also a convergence of moments of orders $p \geq0$.
\end{theorem}

%The convergence (\ref{nblambda}) was already announced by
Note that a result related to (\ref{nblambda}) is announced in \cite{GIM08} (remark following Theorem 3.1
therein).

% in who mention that it could
%probably be proved, under the regular variation assumption
%(\ref{eq:6}) on $\Lambda$, by adapting the proof they propose to
%describe the behavior of the total number of collisions in
%$\Lambda$-coalescents with finite measure $x^{-2}\Lambda(\mathrm dx)$.

%Our approach here is
%completely different and provides moreover the behavior of the entire
%process of numbers of collisions $X_n$.

Of course, the statement of Theorem \ref{theolambda} remains true if we
simultaneously replace $h$ and $\psi$ in~(\ref{eq:6}) and (\ref{eq:7})
with $ch$ and $c\psi$ for any $c>0$. Also, the statement remains true
if we change $h$ with any of its equivalents at $0$ in
(\ref{coalambda}) or (\ref{nblambda}). Theorem \ref{theolambda}
specialises to yield the following results on beta coalescents. Recall
that the beta-coalescent with parameters $a,b>0$, also denoted by
$\beta(a,b)$-coalescent, is the $\Lambda$-coalescent associated with
the measure
\[
\Lambda(\d{x})=\frac{\Gamma(a+b)}{\Gamma(a)\Gamma
(b)}x^{a-1}(1-x)^{b-1}\,\d{x} \,\ind_{[0,1]}(x) .
\]

\begin{corollary} For the beta-coalescent $\beta(a,b)$ with parameters
$1<a<2$ and $b>0,$ the process of numbers of collisions satisfies
\[
\biggl( \frac{X_n(\lfloor n^{2-a}t\rfloor)}{n}, t \geq0\biggr)
\mathop{\stackrel{(d)}{\longrightarrow}}_{n\to\infty}
\exp(-\xi_{\tau}),
\]
where $\xi$ is a subordinator with Laplace exponent
\[
\psi(\lambda)=\frac{2-a}{\Gamma(a)} \int_0^{\infty}
(1-\mathrm{e}^{-\lambda
y})\frac{ \mathrm{e}^{-by}}{(1-\mathrm{e}^{-y})^{3-a}}\,\d{y}
\]
and $\tau$ the time change defined from $\xi$ by (\ref{rho}),
replacing there $\gamma$ with $2-a$. Moreover, the total number $A_n$
of collisions in such a beta-coalescent satisfies, jointly with the
previous convergence,
\[
\frac{A_n}{n^{2-a}}\mathop{\stackrel{(d)}{\longrightarrow}}_{n\to\infty}\int_0^{\infty} \exp\bigl(-(2-a) \xi_r\bigr)\,\d{r} .
\]
The convergence of all positive moments also holds.
\end{corollary}

When $b=2-a$, we know from the particular form of the Laplace exponent
of $\xi$ that the range of $\exp(-\xi)$ is identical in law with the
zero set of a Bessel bridge of dimension $2-2b$ (see \cite{GPreg}).
When, moreover, $b\in(0,1/2]$, the time changed process
$\exp(-\xi_{\tau})$ is distributed as the tagged fragment in a
$1/(1-b)$-stable fragmentation (with a dislocation measure suitably
normalized). More generally, when $b\in(0,1)$ and $a>1+b$, the time
changed process $\exp(-\xi_{\tau})$ is distributed as the tagged
fragment in a Poisson--Dirichlet fragmentation with a dislocation
measure proportional to $\mathit{PD}^*(1-b,a+b-3)$ as defined in \cite{HPW},
Section 3. In such cases, the Laplace exponent of $\xi$ can be
explicitly computed. See Corollary 8 of \cite{HPW}.

When $b=1$ (and still $1<a<2$), the asymptotic behavior of $A_n$ is
proved by Iksanov and M\"{o}hle in \cite{IM08}, using there the
connection with this model and random walks with a barrier. As
mentioned at the end of the previous section, the limit random variable
$\int_0^{\infty}\exp(-(2-a)\xi_t)\,\mathrm{d}t$ is then distributed as
$(a-1)\tau_{2-a}^{a-2}$, where $\tau_{2-a}$ is a $(2-a)$-stable
variable, with Laplace transform $\mathbb E[\exp(-\lambda
\tau_{2-a})]=\exp(-\lambda^{2-a})$.

Besides, Iksanov, M\"{o}hle and co-authors obtain various results on the
asymptotic behavior of $A_n$ for beta coalescents when $a \notin
(1,2)$. See \cite{IMM09} for a summary of these results.

%s2.3 ###
\subsection{Regenerative compositions}\label{sec:numb-blocks-regen}

A \textit{composition} of $n \in\mathbb N$ is a sequence
$(c_1,c_2,\ldots,c_k)$, $c_i \in\mathbb N$, with sum $\sum_{i=1}^k c_i =
n$. The integer $k$ is called the \textit{length} of the composition.
If $X_n$ is a Markov chain taking values in $\mathbb Z_+$,
\textit{strictly decreasing} on $\mathbb N$ and such that $X_n(0) = n$,
the random sequence
\[
C^{(n)}_i:=X_n(i-1)-X_n(i),\qquad 1 \leq i \leq K^{(n)}:=\inf\{k\dvtx X_n(k)=0\},
\]
clearly defines a random composition of $n$, of length $K^{(n)}$.
Thanks to the Markov property of $X$, the random sequence
$(C^{(n)},n\geq1)$ has the following \textit{regenerative} property:
\[
\bigl(C^{(n)}_2,C^{(n)}_3,\ldots,C^{(n)}_{K^{(n)}} \bigr)\quad \mbox{conditional
on}\quad \bigl\{C^{(n)}_1=c_1 \bigr\} \stackrel{\mathrm{law}}= C^{(n-c_1)}\qquad
\forall 1 \leq c_1 <n.
\]
This is called a \textit{regenerative composition}. Conversely,
starting from a regenerative composition $(C^{(n)},n\geq1)$, we build,
for each $n \geq1$, a strictly decreasing Markov chain $X_n$ starting
at $n$ by setting
\[
X_n(k) = n -\sum_{i=1}^k C^{(n)}_i,\qquad 1 \leq k \leq K^{(n)},\quad \mbox{and}\quad
X_n(k)=0 \qquad\mbox{for } k \geq K^{(n)}.
\]
The transition probabilities of the chain are $p_{n,k}=\mathbb
P(C^{(n)}_{1} = n -k)$ for $0 \leq k <n$, $p_{n,n}=0$ for $n \geq1$
and $p_{0,0}=1$.

Regenerative compositions have been studied in great detail by Gnedin
and Pitman \cite{GPreg} under the additional following
\textit{consistency} property: For all $n\geq2$, if $n$ balls are
thrown at random into an ordered series of boxes according to
$C^{(n)}$, then the composition of $n-1$ obtained by deleting one ball
uniformly at random is distributed according to $C^{(n-1)}$. Gnedin
and Pitman \cite{GPreg} show in particular that regenerative consistent
compositions can be constructed via (unkilled) subordinators through
the following procedure. Let $\xi$ be such a subordinator and
$(U_i,i\geq1)$ be an independent sequence of i.i.d. random variables
uniformly distributed on $(0,1)$. Construct from this an ordered
partition of $[n]$, say, $(B^{(n)}_1,\ldots,B^{(n)}_{K^{(n)}})$, by
declaring that $i$ and $j$ are in the same block if and only if $U_i$
and $U_j$ are in the same open interval component of
$[0,1]\backslash\{1-\exp(-\xi_t), t \geq0 \}^{\mathrm{cl}}$. The order
of blocks is naturally induced by the left-to-right order of open
interval components. Then $((\#B^{(n)}_1,\ldots,\#B^{(n)}_{K^{(n)}} ), n
\geq1)$ defines a regenerative consistent composition. Conversely,
each regenerative consistent composition can be constructed in that way
from a subordinator.

In cases where the subordinator has no drift and its L\'{e}vy measure
$\omega$ has a tail that varies regularly at $0$, that is, $\overline
\omega(x):=\int_x^{\infty}\omega(\mathrm dy)=x^{-\gamma}\ell(x),$
where $\gamma\in(0,1)$ and $\ell$ is slowly varying at $0$, Gnedin,
Pitman and Yor \cite{GPYalpha} show that
%begin{theorem}
%
\[
\frac{K^{(n)}}{\Gamma(1-\gamma)n^{\gamma}\ell(1/n)}
\stackrel{\mathrm{a.s.}}\rightarrow\int_0^{\infty}\exp(-\gamma\xi_r)
\,\mathrm{d}r.
\]
The duality between regenerative compositions and strictly decreasing
Markov chains, coupled with Theorem \ref{cvdeathtimes}, allows us to
extend this result by Gnedin, Pitman and Yor to the largest setting of
regenerative compositions that do not necessarily follow the
consistency property, provided hypothesis (H) holds. Note,
however, that in this more general context we can only obtain a
convergence in distribution.

Let us check here that in the consistent cases, the assumption of
regular variation on the tail of the L\'{e}vy measure associated with
the composition entails (H). Following \cite{GPreg}, the
transition probabilities of the associated chain $X$ are then given by
\[
p_{n,k}=\mathbb P\bigl(C^{(n)}_{1} = n-k\bigr)=\frac{1}{Z_n} \pmatrix{n\cr k}
\int_0^{1}x^{k}(1-x)^{n-k} \tilde\omega(\mathrm dx),\qquad 0 \leq k \leq n-1,
\]
where $\tilde\omega$ is the push-forward of $\omega$ by the
mapping $x \mapsto\exp(-x)$ and $Z_n$ is the normalizing constant
$Z_n= \int_0^{1}(1-x^n) \tilde\omega(\mathrm dx). $ It is easy to see
that (\ref{Hyp}) is satisfied with $a_n=Z_n$ and $\mu(\mathrm dx)=(1-x)
\tilde\omega(\mathrm dx)$ since the distributions $(p_{n,k}, 0 \leq k
\leq n-1), n \geq1,$ are mixtures of binomial-type distributions (we
refer to the proof of Proposition~\ref{sec:scaling-limits-non-3} or to
that of the forthcoming Lemma \ref{lemlamb2} for detailed arguments in
a similar context). The Laplace transform defined via $\mu$ by
(\ref{eq:1}) is then that of a subordinator with L\'{e}vy measure
$\omega$, no drift and killing rate $\mathtt{k}=0$. Besides, by
Karamata's Tauberian theorem \cite{BGT}, Theorem 1.7.1$'$, the assumption
$\overline\omega(x) = x^{-\gamma}\ell(x)$, $\gamma\in(0,1),$ where
$\ell$ is slowly varying at $0$, implies that
\[
Z_n= \int_0^{\infty}(1-\mathrm{e}^{-nx}) \omega(\mathrm dx) = n
\int_0^{\infty}\mathrm{e}^{-nu} u^{-\gamma}\ell(u) \,\mathrm{d}u \sim
\Gamma(1-\gamma)n^{\gamma}\ell(1/n)\qquad \mbox{as }n \rightarrow\infty,
\]
and we have indeed (H) with the correct parameters
$(a_n,n\geq1)$ and $\mu$.

Last, we rephrase Theorem \ref{CvChain} in terms of regenerative
compositions.
%To any regenerative composition $(C^{(n)}, n \geq1)$, we associate
%the random probability measure $$ m_{C^{(n)}}:=\sum_{k=1}^{K^{(n)}}

\begin{theorem}
\label{theocomp}
Let $(C^{(n)}, n \geq1)$ be a regenerative
composition.
\begin{longlist}[(ii)]
\item[(i)] Assume that it is consistent, constructed via
a subordinator $\xi$ with no drift and a L\'{e}vy measure with a tail
$\overline\omega$ that varies regularly at $0$ with index $-\gamma$,
$\gamma\in(0,1)$. Then,
\[
\biggl( \sum_{k \leq t \overline\omega(1/n) \Gamma(1-\gamma)}
\frac{C^{(n)}_k}{n}, t \geq0 \biggr) \stackrel{ \mathrm{(d)}}
\rightarrow\bigl(1-\exp\bigl(-\xi_{\tau(t)}\bigr),t\geq0 \bigr),
\]
%
%$$ ( m_{C^{(n)}}([0,t]), t \geq0 ) \overset{
%their cumulative distribution functions that converge in
%distribution, for the Skorokhod topology, to
%$(1-\exp(-\xi_{\tau(t)}),t\geq0)$,
where $\tau$ is the usual time change defined as the inverse of $t
\mapsto\int_0^t \exp(-\gamma\xi_r) \,\mathrm{d}r$.
\item[(ii)] When the regenerative composition is
non-consistent, assume that $\mathbb E[C_1^{(n)}]/n$ varies regularly
as $n \rightarrow\infty$ with index $-\gamma$, $\gamma\in(0,1]$ and
that
\[
\frac{ \mathbb E[ C_1^{(n)}f(1-C_1^{(n)}/n )]} {\mathbb E[C_1^{(n)}]}
\rightarrow\int_{[0,1]}f(x)\mu(\mathrm dx)
\]
for a probability measure $\mu$ on $[0,1]$ and all continuous functions
$f \dvtx [0,1] \rightarrow\mathbb R_+$. Then,
\[
\biggl( \sum_{k \leq t n/ \mathbb E[C_1^{(n)}]} \frac{C^{(n)}_k}{n}, t
\geq0 \biggr) \stackrel{ \mathrm{(d)}} \rightarrow
\bigl(1-\exp\bigl(-\xi_{\tau(t)}\bigr),t\geq0 \bigr),
\]
where $\xi$ is the subordinator with Laplace exponent defined via
$\mu$ by (\ref{eq:1}) and $\tau$ the usual time change.
\end{longlist}
\end{theorem}

As was pointed out to us by a referee, the assertion (i) in this
statement actually holds in the almost-sure sense. This is an easy
consequence of \cite{GPYalpha}, Theorem 4.1.

%%%%%%%%%%%%%%%%%%
%s3 ###
\section{Preliminaries}
\label{prem}
%%%%%%%%%%%%%%%%%%

Our goal now is to prove Theorems \ref{CvChain} and \ref{cvdeathtimes}
and Proposition~\ref{sec:scaling-limits-non-2}. We start in this
section with some preliminaries. From now on and until the end of
Section \ref{Markovchain}, we suppose that assumption (H) is in
force. Consider the generating function defined for all $\lambda\geq
0$ by
%
%e13 ###
\begin{equation}
\label{fonge}
G_n(\lambda)=\sum_{k=0}^n
\biggl(\frac{k}{n}\biggr)^{\lambda}p_{n,k}=\mathbb E \biggl[ \biggl(\frac{X_n(1)}{n}
\biggr)^{\lambda}\biggr]
\end{equation}
with the convention $G_0(\lambda)=0$. Then
%
%e14 ###
\begin{equation}\label{eq:11}
1-G_n(\lambda)=\sum_{k=0}^n\biggl(1- \biggl(\frac{k}{n}\biggr)^{\lambda}\biggr)p_{n,k}
=\int_{[0,1]}[\lambda]_x(1-x) p^*_{n}(\d{x}) ,
\end{equation}
where $[\lambda]_x$ was defined around (\ref{eq:3}). Thanks to (H), we immediately get
%
%e15 ###
\begin{equation}\label{eq:4}
a_n\bigl(1-G_n(\lambda)\bigr) \mathop{\longrightarrow}\limits_{n\to\infty}\psi
(\lambda)
,\qquad \lambda> 0 ,
\end{equation}
the limit being the Laplace exponent defined at (\ref{eq:1}). In fact,
if this convergence holds for every $\lambda>0$, then (\ref{Hyp})
holds.

\begin{proposition}\label{sec:preliminaries}
Assume that there exists a sequence of the form $a_n =
n^{\gamma}\ell(n),n\geq1$ for some slowly varying function
$\ell\dvtx \mathbb R^+ \rightarrow(0,\infty)$, such that (\ref{eq:4}) holds
for some function $\psi$ and every $\lambda> 0$, or only for an
infinite set of values of $\lambda\in(0,\infty)$ having at least one
accumulation point. Then there exists a unique finite measure $\mu$ on
$[0,1]$ such that $\psi(\lambda)=\int_0^1[\lambda]_x\mu(\d{x})$ for
every $\lambda>0$, and \textup{(H)} holds for the sequence $(a_n,n\geq1)$
and the measure $\mu$.
\end{proposition}

\begin{pf}
For any given $\lambda>0$, the function $x
\mapsto[\lambda]_x$ is bounded from below on $ [0,1] $ by a positive
constant $c_\lambda>0$. Therefore, if (\ref{eq:4}) holds, then, using
(\ref{eq:11}), we obtain that
\[
\sup_{n\geq1}\int_{[0,1]}a_n(1-x)p^*_n(\d{x})\leq
\frac{1}{c_\lambda}\sup_{n\geq1}a_n\bigl(1-G_n(\lambda)\bigr)<\infty.
\]
Together with the fact that the measures $a_n(1-x)p^*_n(\d{x}),n\geq1$
are all supported on $[0,1]$, this implies that all subsequences of
$(a_n(1-x)p^*_n(\d{x}),n\geq1)$ have a weakly convergent subsequence.
Using (\ref{eq:4}) again, we see that any possible weak limit $\mu$
satisfies $\psi(\lambda)=\int_0^1[\lambda]_x\mu(\d{x})$. This function
is analytic in $\lambda>0$, and uniquely characterizes $\mu$. The same
holds if we only know this function on an infinite subset of
$(0,\infty)$ having an accumulation point, by analytic continuation.
\end{pf}

For some technical reasons, we need for the proofs to work with
sequences $(a_n,n\geq0)$ rather than sequences indexed by $\mathbb
N$. We therefore complete all the sequences $(a_n,n\geq1)$ involved in
(H) or (\ref{eq:4}) with an initial term $a_0=1$. This is
implicit in the whole Sections \ref{prem} and \ref{Markovchain}.

%%%%%%%%%%%%%%%%%%%%%%%%%%%%%%%%
%s3.1 ###
\subsection{Basic inequalities}\label{sec:crucial-inequalities}
%%%%%%%%%%%%%%%%%%%%%%%%%%%%%%%%

Let $\lambda>0$ be fixed.
% is fixed and $c_i(\lambda)$, $i=1,2,3,4$ are
%strictly positive finite reals, that depend in general on $\lambda$,
%but not on $n$.
By (\ref{eq:4}), there exists
%together with the fact that $a_n>0$ for every
%$n\geq0$, entails the existence of
a finite constant $c_1(\lambda)>0$ such that for every $n\geq0$,
%
%e16 ###
\begin{equation}
\label{crucialinequalities}
1-G_n(\lambda) \leq\frac{c_1(\lambda)}
{a_n} .
\end{equation}
In particular, $G_n(\lambda)>1/2$ for $n$ large enough. Together with
the fact that $a_n>0$ for every $n\geq0$, this entails the existence
of an integer $n_0(\lambda)\geq0$ and finite constants
$c_2(\lambda),c_3(\lambda)>0$ such that, for every $n\geq
n_0(\lambda)$,
%
%e17 ###
\begin{equation}
\label{ceq4}
-\ln(G_n(\lambda) ) \leq\frac{c_2(\lambda)} {a_n}
\leq
c_3(\lambda) .
\end{equation}
When, moreover, $p_{n,n}<1$ for all $n \geq1$ (or, equivalently,
$G_n(\lambda)<1$ for all $n \geq1$), we obtain the existence of a
finite constant $c_4(\lambda)>0$ such that, for every $n\geq1$
%
%e18 ###
\begin{equation}
\label{ceq3}
\ln( G_{n}(\lambda) ) \leq-\frac{c_4(\lambda)}{a_n} .
\end{equation}
Last, since $(a_n,n\geq0)$ is regularly varying with index $\gamma$
and since $a_n>0$ for all $n \geq0$, we get from Potter's bounds
\cite{BGT}, Theorem 1.5.6, that for all $\varepsilon>0$, there exist
finite positive constants $c'_1(\varepsilon)$ and $c'_2(\varepsilon)$
such that, for all $1 \leq k \leq n$
%
%e19 ###
\begin{equation}
\label{ceq2}
c'_1(\varepsilon) \biggl( \frac{n}{k} \biggr)^{\gamma-\varepsilon}
\leq\frac{a_n}{a_k} \leq c'_2(\varepsilon) \biggl( \frac{n}{k}\biggr)^{\gamma+\varepsilon}.
\end{equation}

%%%%%%%%%%%%%%%
%s3.2 ###
\subsection{Time changes}
\label{tc}
%%%%%%%%%%%%%%%
Let $f\dvtx [0,\infty) \rightarrow[0,1]$ be a c\`{a}dl\`{a}g non-increasing
function. We let $\sigma_f:=\inf\{t \geq0 \dvtx f(t)=0 \}$, with the
convention $\inf\{\varnothing\}=\infty$.
% We will write $\sigma$ instead of
%$\sigma_f$ if there is no ambiguity.
Now fix $\gamma>0$. For $0\leq t<\sigma_f$, we let
\[
\tau_f(t):=\int_0^t f(r)^{-\gamma}\,\d{r} ,
\]
and $\tau_f(t)=\infty$ for $t\geq\sigma_f$. Then $(\tau_f(t),t\geq0)$
is a right-continuous, non-decreasing process with values in
$[0,\infty]$, and which is continuous and strictly increasing on
$[0,\sigma_f)$. Note that $\tau_f(\sigma_f-)=\int_0^{\sigma_f}
f(r)^{-\gamma}\,\d{r}$ might be finite or infinite.
We set
\[
\tau_f^{-1}(t)=\inf\{u\geq0\dvtx \tau_f(u)>t\} ,\qquad t\geq0 ,
\]
which
defines a continuous, non-decreasing function on $\mathbb R_+$, that is
strictly increasing on $[0,\tau_f(\sigma_f-))$, constant equal to
$\sigma_f$ on $[\tau_f(\sigma_f -),\infty)$, with limit
$\tau^{-1}_f(\infty)=\sigma_f$. The functions $\tau_f$ and
$\tau_f^{-1}$, respectively restricted to $[0,\sigma_f)$ and
$[0,\tau_f(\sigma_f-))$, are inverses of each other. The function
$\tau_f$ is recovered from $\tau_f^{-1}$ by the analogous formula
$\tau_f(t)=\inf\{u \geq0\dvtx \tau_f^{-1}(u) > t\}$, for any $ t \geq0$.

We now consider the function
\[
g(t):=f(\tau_f^{-1}(t) ),\qquad t \geq0,
\]
which is also c\`{a}dl\`{a}g, non-increasing, with values in $[0,1]$, and
satisfies $\sigma_g=\tau_f(\sigma_f-)$. Note also that $f(t)=g
(\tau_f(t) ) ,$ $t \geq0 , $ where, by convention, $g(\infty)=0$.
Finally, we have
\[
\d\tau_f(t)=f(t)^{-\gamma} \,\d{t}\qquad \mbox{on } [0,\sigma_f)
\]
and
\[
\tau_f^{-1}(t)=f(\tau_f^{-1}(t))^{\gamma}\,\d{t}=g(t)^{\gamma}\,\d{t}
\qquad\mbox{on }[0,\sigma_g).
\]
Now, for $c>0$, we will often
use the change of variables $u=\tau_f(r/c)$ to get that when
$g(t)>0$ (i.e., $t<\sigma_g$), for any measurable, non-negative
function $h$,
%
%e20 ###
\begin{equation}
\label{chvar} \int_0^{c\tau_f^{-1}(t)} h(f(r/c) ) \,\mathrm{d}r=c \int_0^t
h(g(u) )g(u)^{\gamma}\,\mathrm{d}u.
\end{equation}
In particular $ \tau_f^{-1}(t)=\int_0^{t}g(r)^{\gamma} \,\mathrm{d}r $ for
$t<\tau_f(\sigma_f-)$. This remains true for $t\geq\tau_f(\sigma_f-)$
since $g(t)=0$ for $t \geq\tau_f(\sigma_f-)$. Consequently, $
\tau_f^{-1}(t)=\int_0^{t}g^{\gamma}(r) \,\mathrm{d}r $ for all $t \geq0$
and
%
%e21 ###
\begin{equation}
\label{sigma} \sigma_f=\int_0^{\infty} g(r)^{\gamma} \,\d{r}.
\end{equation}
This also implies that $\tau_f(t)=\inf\{u \geq0\dvtx \int_0^{u}g^{\gamma}(r) \,\mathrm dr> t \}$ for every $t \geq0.$

%and we can reverse the roles of $f$ and $g$.

%%%%%%%%%%%%%%%%%%
%s3.3 ###
\subsection{Martingales associated with $X_n$}
\label{martingales}
%%%%%%%%%%%%%%%%%%

%There is a very classical way of associating martingales with any
%Markov chain. Two particular forms are going to be useful for our
%purposes.

We finally recall the very classical fact that if $P$ is the
transition function of a Markov chain $X$ with countable state space
$M$, then for any non-negative function $f$, the process defined by
\[
f(X(k))+\sum_{i=0}^{k-1} (\mathrm{Id}-P)f(X(i)) ,\qquad k\geq0,
\]
is a martingale, provided all the terms of this process are integrable.
When, moreover, $f^{-1}(\{ 0\})$ is an absorbing set (i.e., $f(X(k))=0$
implies $f(X(k+1))=0$), the process defined by
\[
f(X(k))\prod_{i=0}^{k-1} \frac{f(X(i))}{Pf(X(i))} ,\qquad k\geq0 ,
\]
with the convention $0\cdot\infty=0$ is also a martingale
(absorbed at $0$), provided all the terms are integrable. From this, we
immediately obtain the following.

\begin{proposition}
\label{propmartingales}
For every $\lambda>0$ and every integer $n
\geq 1$, the processes defined by
%
%e22 ###
\begin{equation}\label{eq:13}
\biggl(\frac{X_n(k)}{n} \biggr)^{\lambda}+\sum_{i=0}^{k-1} \biggl(\frac{X_n(i)}{n}
\biggr)^{\lambda}\bigl(1-G_{X_n(i)}(\lambda)\bigr) ,\qquad k\geq0,
\end{equation}
and
%
%e23 ###
\begin{equation}\label{eq:14}
\Upsilon^{(\lambda)}_n(k)= \biggl(\frac{X_n(k)}{n} \biggr)^{\lambda} \Biggl(
\prod_{i=0}^{k-1}G_{X_n(i)}(\lambda)\Biggr)^{-1} ,\qquad k\geq0,
\end{equation}
are martingales with respect to the filtration generated by $X_n$, with
the convention that $\Upsilon^{(\lambda)}_n(k)=0$ whenever $X_n(k)=0$.
\end{proposition}

%To justify the convention adopted for $\Upsilon^{(2)}$, note that if
%$K$ is the first time where $X_n$ attains $0$, then $X_n(i)>0$ for
%$i<K$ so that $G_{X_n(i)}(\lambda)> 0$. Therefore, $0$ is an absorbing
%state for $$

%%%%%%%%%%%%%%%%%%%%%%%%%%%%
%s4 ###
\section{Scaling limits of non-increasing Markov chains}
\label{Markovchain}
%%%%%%%%%%%%%%%%%%%%%%%%%%%%

We now start the proof of Theorems \ref{CvChain} and \ref{cvdeathtimes}
and Proposition~\ref{sec:scaling-limits-non-2}. As mentioned before,
this is done by first establishing tightness for the processes
$Y_n(t)=X_n(\lfloor a_n t\rfloor)/n,t\geq0$. We recall that
(H) is assumed throughout the section, except in the last
subsection, which is devoted to the proof of Proposition~\ref{sec:scaling-limits-non-3}.

%%%%%%%%%%%%%%%
%s4.1 ###
\subsection{Tightness}
\label{Tightness}
%%%%%%%%%%%%%%%

\begin{lemma}
\label{LemTightness}
The sequence $(Y_n,n\geq0)$ is tight with respect
to the Skorokhod topology.
\end{lemma}

Our proof is based on Aldous' tightness criterion, which we first
recall.

\begin{lemma}[(Aldous' tightness criterion
\cite{billingsley99}, Theorem
16.10)]
\label{Aldous}
Let $(F_n, n \geq0)$ be a sequence of
$\mathcal{D}$-valued stochastic processes and for all $n$ denote by
$\mathcal{J}(F_n)$ the set of stopping times with respect to the
filtration generated by $F_n$. Suppose that for all fixed $t>0$,
$\varepsilon>0$
\begin{eqnarray*}
\label{Aldouscriterion}
&&\phantom{i}\mathrm{(i)}\quad \lim_{a \rightarrow\infty}\limsup_{n \rightarrow\infty} \mathbb P\Bigl( \sup_{s \in[0,t]} F_n(s)>
a\Bigr)=0;
\\
&&\mathrm{(ii)}\quad \lim_{\theta_0 \rightarrow0} \limsup_{n
\rightarrow\infty} \sup_{T \in\mathcal{J}(F_n), T \leq t} \sup_{0\leq\theta\leq\theta_0} \mathbb P\bigl(|F_n(T)-F_n(T+\theta)| >\varepsilon\bigr)=0,
\end{eqnarray*}
then the sequence $(F_n,n\geq0)$ is tight with respect to the Skorokhod
topology.
\end{lemma}

\begin{pf*}{Proof of Lemma \ref{LemTightness}}
Part (i) of
Aldous' tightness criterion is obvious since $Y_n(t) \in[0,1]$, for
every $n \in\mathbb Z_+, t \geq0$. To check part (ii), consider some
$\lambda>\max(\gamma,1)$, where $\gamma$ denotes the index of regular
variation of $(a_n,n\geq0)$.
%With a slight abuse of notations, we denote by $T_n$ the set of
%stopping times with respect to the filtration generated by $Y_n$.
Then, on the one hand, for all $n \geq1$, since the process $Y_n$ is
non-increasing and $\lambda\geq1$, we have for all (possibly random)
times $T$ and all $\theta\geq0$,
\[
|Y_n(T)-Y_n(T+\theta)|^{\lambda} \leq
Y_n^{\lambda}(T)-Y_n^{\lambda}(T+\theta).
\]
On the other hand, let $T$ be a bounded stopping time in
$\mathcal{J}(Y_n)$. Then $\lfloor a_nT\rfloor$ is a stopping time with
respect to the filtration generated by $X_n$. Applying Doob's optional
stopping theorem to the martingale (\ref{eq:13}) yields, for every
$\theta\geq0$,
\begin{eqnarray*}
\mathbb E[Y_n^{\lambda}(T)-Y_n^{\lambda}(T+\theta)]
% \rfloor)-X_n^{\lambda}(\lfloor a_n(T+\theta)\rfloor)]
&=& n^{-\lambda} \mathbb E\Biggl[ \sum_{i=\lfloor a_nT\rfloor}^{\lfloor
a_n(T+\theta)\rfloor-1} X_n^{\lambda}(i) \bigl(1-G_{X_n(i)}(\lambda) \bigr) \Biggr]
\\
& \leq& c_1(\lambda) n^{-\lambda} \mathbb E\Biggl[ \sum_{i=\lfloor
a_nT\rfloor}^{\lfloor a_n(T+\theta)\rfloor-1} \frac{X_n^{\lambda}(i)}
{a_{X_n(i)}} \Biggr],
\end{eqnarray*}
where we used (\ref{crucialinequalities}) at the last step. Next,
since $\lambda>\gamma$, $X_n^{\lambda}(i) / a_{X_n(i)} \leq
c'_2(\lambda-\gamma) n^{\lambda}/a_n$ for all $n,i \geq0$, where
$c'_2(\varepsilon)$ was introduced in (\ref{ceq2}) (note that the
inequality is obvious when $X_n(i)=0$, since $a_0>0$). Hence, for every
bounded $T \in\mathcal{J}(Y_n)$ and $\theta\geq0$,
\begin{eqnarray*}
\mathbb E[|Y_n(T)-Y_n(T+\theta)|^{\lambda}] &\leq& \mathbb
E[Y_n^{\lambda}(T)-Y_n^{\lambda}(T+\theta)] \\ & \leq&
\frac{c_1(\lambda) c'_2(\lambda-\gamma) }{a_n} \mathbb E[ \lfloor
a_n(T+\theta)\rfloor-\lfloor a_nT\rfloor] \\ &\leq& c_1(\lambda)
c'_2(\lambda-\gamma) (\theta+a_n^{-1} ),
\end{eqnarray*}
which immediately yields (ii) in Aldous' tightness criterion.
\end{pf*}

%%%%%%%%%%%%%%%%%%%%%
%s4.2 ###
\subsection{Identification of the limit}
\label{Limit}
%%%%%%%%%%%%%%%%%%%%%

We now want to prove uniqueness of the possible limits in distribution
of subsequences of $Y_n, n \geq0$.
%Thanks to Lemma \ref{Tightness},
Let $(n_k,k\geq0)$ be a strictly increasing sequence, such that the
process $Y_n$ converges in distribution to a limit $Y'$ when $n$ varies
along $(n_k)$. To identify the distribution of $Y'$, recall the
definition of $Z_n=(Y_n(\tau_n^{-1}(t)),t\geq0)$ at the end of Section~\ref{secscal}. From the discussion in Section~\ref{tc}, we have
\[
Z_n(t)=Y_n \biggl(\int_0^{t} Z_n(r)^{\gamma} \,\mathrm dr\biggr) ,\qquad t \geq0 .
\]
As in Section \ref{tc}, let
$\tau_{Y'}(u)=\int_0^u Y'(r)^{-\gamma}\,\d{r}$ if $Y'(u)>0$ and
$\tau_{Y'}(u)=\infty$ otherwise, and let $Z'$ be the process defined by
\[
Z'(t)=Y'(\tau_{Y'}^{-1}(t)) ,
\]
where $\tau_{Y'}^{-1}(t)=\inf\{u\geq
0: \tau_{Y'}(u)>t\}$, so that
\[
Z'(t)=Y'\biggl(\int_0^{t} Z'(r)^{\gamma} \,\mathrm{d}r\biggr) ,\qquad t\geq0 .
\]
Then, as a consequence of \cite{EK}, Theorem 1.5,
Chapter 6 (it is in fact a consequence of a step in the proof of this
theorem rather than its exact statement), the convergence in
distribution of $Y_n$ to $Y'$ along $(n_k)$ entails that of $(Y_n,Z_n)$
to $(Y',Z')$ in $\mathcal{D}^2$ along the same subsequence, provided
the following holds:
\[
\sigma_{Y'}=\inf\{ s \geq0 \dvtx Y'(s)^{\gamma}=0 \}=\inf\biggl\{ s \geq0 \dvtx
\int_0^s Y'(u)^{-\gamma}\,\mathrm d u = \infty
\biggr\}=\lim_{\varepsilon\rightarrow0} \inf\{ s \geq0 \dvtx Y'(s)^{\gamma} <
\varepsilon\},
\]
which is obviously true here since $Y'$ is a.s. c\`{a}dl\`{a}g non-increasing.
Therefore, the proof of Theorem \ref{CvChain} and Proposition~\ref{sec:scaling-limits-non-2} will be completed provided we show the
following.

\begin{lemma}
\label{lemmachar} The process $Z'$ has same distribution as
$Z=(\exp(-\xi_t),t\geq0)$, where $\xi$ is a subordinator with Laplace
exponent $\psi$.
\end{lemma}

To see that this entails Proposition~\ref{sec:scaling-limits-non-2}
(hence Theorem \ref{CvChain}), note that $\tau_{Y'}(t)=\inf\{u \geq
0 \dvtx \int_0^{u} Z'(r)^{\gamma} \,\mathrm dr>t\}$ and $Y'(t)=Z'(\tau_{Y'}(t))$,
for $t\geq0$, as detailed in Section \ref{tc}. So the previous lemma
entails that the only possible limiting distribution for $(Y_n)$ along
a subsequence is that of $Y$ as defined in (\ref{eq:2}). Since
$(Y_n,n\geq0)$ is a tight sequence, this shows that it converges in
distribution to $Y$, and then that $(Y_n,Z_n)$ converges to $(Y,Z)$,
entailing Proposition~\ref{sec:scaling-limits-non-2}.

To prove Lemma \ref{lemmachar}, we need a pair of results on Skorokhod
convergence, which are elementary and left to the reader. The first
lemma is an obvious consequence of the definition of Skorokhod
convergence. The second one can be proved, for example, by using
Proposition~6.5 in \cite{EK}, Chapter~3.

\begin{lemma}
\label{lemmaSko1}
Suppose that $f_n \rightarrow f$ on $\mathcal{D} $
and that $(g_n,n\geq0)$ is a sequence of c\`{a}dl\`{a}g non-negative functions
on $[0,\infty)$ converging uniformly on compacts to a continuous
function
$g$. Then $f_ng_n \rightarrow fg$ on $\mathcal{D}$.
\end{lemma}

\begin{lemma}
\label{lemmaSko2}
Suppose that $f_n, f$ are non-increasing,
non-negative functions in $\mathcal{D}$ such that $f_n \rightarrow f$.
Let $\varepsilon>0$ be such that there is at most one $x \in
[0,\infty)$ such that
$f(x)=\varepsilon$. Define
\[
t_{n,\varepsilon}:=\inf\{t \geq0\dvtx f_n(t) \leq\varepsilon\} \quad\mbox{and}\quad
t_{\varepsilon}:=\inf\{t \geq0\dvtx f(t) \leq\varepsilon\}
\]
(which can be infinite).
Then it holds that $ t_{n,\varepsilon} \rightarrow t_{\varepsilon}$ as
$n\to\infty$, and if $f(t_{\varepsilon}-)>\varepsilon$ or
$f(t_{\varepsilon}-)=f(t_{\varepsilon})$, then
\[
\bigl(f_n(t \wedge t_{n,\varepsilon}),t\geq0 \bigr) \rightarrow\bigl(f(t \wedge
t_{\varepsilon}),t\geq0 \bigr).
\]
\end{lemma}

\begin{pf*}{Proof of Lemma \ref{lemmachar}}
Fix $\lambda>0$ and
consider the martingale $(\Upsilon_n^{(\lambda)}(k),k\geq0)$ of
Proposition~\ref{propmartingales}. This is a martingale with respect to
the filtration generated by $X_n$. Therefore, the process
$(\Upsilon_n^{(\lambda)}(\lfloor a_n t\rfloor),t\geq0)$ is a
continuous-time martingale with respect to the filtration generated by~$Y_n$. Next, note that for all $t \geq0$, $\tau_n^{-1}(t)$ is a
stopping time with respect to this filtration, which is bounded (by
$t$). Hence, by Doob's optional stopping theorem, the
process\looseness=-1
%
%e24 ###
\begin{equation}
\label{martingaleMn}
M_n^{(\lambda)}(t)= %
% \{ \begin{array}{l}
%1 \text{ if }\tau_n^{-1}(t)a_n<1 \\
Z_n(t)^{ \lambda} \Biggl( \prod_{i=0}^{\lfloor a_n\tau_n^{-1}(t) \rfloor-1}
G_{X_n(i)}(\lambda)\Biggr)^{-1} ,\qquad t\geq0
% \text{
% if }Z_n(t)> 0\\ 0 \text{ if }Z_n(t)=0 \end{array} .
\end{equation}
(with the usual convention $0\cdot\infty=0$) is a continuous-time
martingale with respect to the filtration generated by $Z_n$.\vadjust{\goodbreak}

We want to exploit the sequences of martingales $(M_n^{(\lambda
)},n\geq0)$ in order to prove that the processes
$(Z'(t)^{\lambda}\exp(\psi(\lambda) t),t\geq0 )$ are (c\`{a}dl\`{a}g)
martingales with respect to the filtration that they generate, for
every $\lambda>0$. It is then easy to check that $-\ln(Z')$ is a
subordinator starting from 0 with Laplace exponent $\psi$.

Using the Skorokhod representation theorem, we may assume that the
convergence of $Z_n$ to $Z'$ along $(n_k)$ is almost sure. We consider
stopped versions of the martingale $M_n^{(\lambda)}$. For all
$\varepsilon>0$ and all $n \geq1$, let
\[
T_{n,\varepsilon}:=\inf\{t \geq0\dvtx Z_n(t) \leq\varepsilon\}\quad \mbox{and}\quad
T_{\varepsilon}:=\inf\{t \geq0\dvtx Z'(t) \leq\varepsilon\},
\]
(which are possibly infinite) and note that $T_{n,\varepsilon}$
(resp. $T_{\varepsilon}$) is a stopping time with respect to the
filtration generated by $Z_n$ (resp. $Z'$).

Let $C_1$ be the set of positive real numbers $\varepsilon>0$ such that
\[
\P\bigl(\exists t_1,t_2\geq0\dvtx t_1\neq t_2, Z'(t_1)=Z'(t_2)=\varepsilon\bigr)>0 .
\]
We claim that this set is
at most countable. Indeed, fix an $\eps>0$ and an integer $K>0$, and
consider the set
\[
B_{\eps,K}=\{\exists t_1,t_2\in[0,
K]\dvtx |t_1-t_2|>K^{-1},Z'(t_1)=Z'(t_2)=\eps\} .
\]
Let $C_{1,K}$ be
the set of numbers $\eps$ such that $\P(B_{\eps,K})>K^{-1}$. If this
set contained an infinite sequence $(\eps_i,i\geq0)$, then by the
reverse Fatou lemma, we would obtain that the probability that
infinitely many of the events $(B_{\eps_i,K},i\geq0)$ occur is at
least $K^{-1}$. Clearly, this is impossible. Therefore, $C_{1,K}$ is
finite for every integer $K>0$. Since $C_1$ is the increasing union
\[
C_1=\bigcup_{K\in\N} C_{1,K} ,
\]
we conclude that it is at most countable. For similar reasons, the set
$C_2$ of real numbers $\varepsilon>0$ such that $\P(
Z'(T_{\varepsilon}-)=\varepsilon>Z'(T_{\varepsilon}) )>0$ is at most
countable.

%(it is implicit that all the assertions here and in the following
%hold almost surely.)

In the rest of this proof, although all the statements and convergences
are in the almost sure sense, we omit the ``a.s.'' in order to have a
lighter presentation. Our goal is to check that for all $\lambda>0$
and all $\varepsilon\notin C_1 \cup C_2$,
\begin{enumerate}[(a)]
\item[(a)] as $n \rightarrow\infty$, the sequence of martingales
$(M_n^{(\lambda)}(t \wedge T_{n,\varepsilon} ), t \geq0 )$
converges to the process $(Z'(t \wedge
T_{\varepsilon})^\lambda\exp(\psi(\lambda) (t \wedge
T_{\varepsilon}
)), t \geq0),$
\item[(b)] the process $(Z'(t \wedge
T_{\varepsilon})^{\lambda}\exp(\psi(\lambda) (t \wedge
T_{\varepsilon} )), t \geq0)$ is a martingale with respect to its
natural filtration,
\item[(c)] the process $(Z'(t)^{\lambda}\exp(\psi(\lambda) t )), t
\geq
0)$ is a martingale with respect to its natural
filtration.
\end{enumerate}

We start with the proof of (a). Fix $\lambda>0$, a positive
$\varepsilon\notin C_1 \cup C_2$ and recall the definition of
$n_0(\lambda)$ in (\ref{ceq4}). Let $n \geq n_0(\lambda)/\varepsilon$.
When $Z_n(t \wedge T_{n,\varepsilon}) > 0$, we can rewrite
%
%e25 ###
\begin{equation}
\label{Mnintegral}
M_n^{(\lambda)}(t \wedge T_{n,\varepsilon}) = \bigl(
Z_n(t \wedge T_{n,\varepsilon})\bigr)^{ \lambda} \exp\biggl(\int_0^{\lfloor
a_n\tau_n^{-1}(t \wedge T_{n,\varepsilon})\rfloor} -\ln\bigl(
G_{X_n(\lfloor
r\rfloor)}(\lambda) \bigr)\,\d{r} \biggr).
\end{equation}
This identity is still true when $Z_n(t \wedge T_{n,\varepsilon}) =0$.
Indeed, even when $Z_n(t \wedge T_{n,\varepsilon}) =0$, for $r <
\lfloor a_n\tau_n^{-1}(t \wedge T_{n,\varepsilon})\rfloor$, it
holds that ${X_n(\lfloor r\rfloor)}\geq n \varepsilon\geq
n_0(\lambda)$. Therefore, by (\ref{ceq4}), the integral involved in
(\ref{Mnintegral}) is well defined and finite. Hence (\ref{Mnintegral})
is valid for all $t \geq0$.

More precisely, as soon as $r < \lfloor a_n\tau_n^{-1}(t \wedge
T_{n,\varepsilon})\rfloor$, we have by (\ref{ceq4}) that
\[
-\ln G_{X_n(\lfloor r\rfloor)}(\lambda) \leq\frac{c_2(\lambda)}
{a_{X_n(\lfloor r\rfloor)}},
\]
which, together with the change of variable identity (\ref{chvar}),
implies that
\[
\int_0^{\lfloor a_n\tau_n^{-1}(t \wedge T_{n,\varepsilon}) \rfloor}
-\ln\bigl( G_{X_n(\lfloor r\rfloor)}(\lambda) \bigr)\,\d{r} \leq
c_2(\lambda)\int_0^{t \wedge T_{n,\varepsilon}} \frac{a_n}
{a_{nZ_n(r)}} Z_n(r)^{\gamma}\,\d r.
\]
Potter's bounds (\ref{ceq2}) and the fact that $Z_n(r)>\varepsilon$ for
$r<t \wedge T_{n,\varepsilon}$ lead to the existence of a finite
constant $c_{\lambda,\varepsilon}$ such that for every $ r<t \wedge
T_{n,\varepsilon}$,
\[
\frac{c_2(\lambda)a_nZ_n(r)^{\gamma}}{a_{nZ_n(r)}} \leq
c_{\lambda,\varepsilon} .
\]
Therefore, for every $t\geq0$,
\[
\int_0^{\lfloor a_n\tau_n^{-1}(t \wedge T_{n,\varepsilon})\rfloor}
-\ln\bigl( G_{X_n(\lfloor r\rfloor)}(\lambda) \bigr)\,\d{r} \leq
c_{\lambda,\varepsilon} (t \wedge T_{n,\varepsilon}) \leq
c_{\lambda,\varepsilon}t .
\]
In particular,
%
%e26 ###
\begin{equation}
\label{MajorationMartingale} M_n^{(\lambda)}(t \wedge
T_{n,\varepsilon}
) \leq\exp(c_{\lambda,\varepsilon}t ) \qquad \forall t \geq0.
\end{equation}
Now we let $n \rightarrow\infty$. Since $\varepsilon\notin C_1 \cup
C_2$, we have, by Lemma \ref{lemmaSko2}, with probability 1,
\[
T_{n,\varepsilon}\rightarrow T_{\varepsilon}\quad \mbox{and}\quad \bigl(Z_n(t
\wedge T_{n,\varepsilon}),t\geq0 \bigr) \rightarrow\bigl(Z'(t \wedge
T_{\varepsilon}),t\geq0 \bigr).
\]
%
%as $n \rightarrow\infty$.
Using (\ref{Mnintegral}) and
%the fact that $n \geq n_0(\lambda)/\varepsilon$ and
Lemma \ref{lemmaSko1}, we see that it is sufficient to prove that
%
%e27 ###
\begin{equation}
\label{CVU}
\biggl( \exp\biggl(\int_0^{\lfloor a_n\tau_n^{-1}(t \wedge
T_{n,\varepsilon})\rfloor} -\ln\bigl( G_{X_n(\lfloor r\rfloor)}(\lambda)
\bigr)\,\mathrm dr \biggr), t \geq0 \biggr) \mathop{\rightarrow}\limits_{n \rightarrow\infty}\bigl(
\exp\bigl(\psi(\lambda)(t \wedge T_{\varepsilon})\bigr), t \geq0 \bigr)
\end{equation}
uniformly on compacts to get the convergence of martingales stated in
(a).

Since we are dealing with non-decreasing processes and the limit is
continuous, it is sufficient to check the pointwise convergence by
Dini's theorem. Fix $t \geq0$. It is well known (see \cite{EK},
Proposition~5.2, Chapter 3) that the Skorokhod convergence implies that
$Z_n(r) \rightarrow Z'(r)$ for all $r$ that is not a jump time of $Z'$,
hence for a.e. $r$. For such an $r$, if $Z'(r) > 0$, we have
$nZ_n(r)\rightarrow\infty$. Hence if $r < t\wedge T_{n,\varepsilon}$,
we have $Z_n(r) \geq\varepsilon$, so that
\[
-\ln\bigl(G_{nZ_n(r)}(\lambda)\bigr) a_n Z_n(r)^{\gamma}
\mathop{\sim}\limits_{n\rightarrow\infty}
\frac{a_n}{a_{nZ_n(r)}}Z_n(r)^{\gamma}\psi(\lambda)
\mathop{\rightarrow}\limits_{n\rightarrow\infty} \psi(\lambda),
\]
using the uniform convergence theorem for slowly varying functions
(\cite{BGT}, Theorem 1.2.1). Moreover, as explained above, the
left-hand side of this expression is bounded from above by
$c_{\lambda,\varepsilon}$ as soon as $n \geq n_0(\lambda
)/\varepsilon$.
This implies, using (\ref{chvar}), that
\[
\int_0^{a_n\tau_n^{-1}(t \wedge T_{n,\varepsilon})} -\ln\bigl(
G_{X_n(\lfloor r\rfloor)}(\lambda) \bigr)\,\d{r}=\int_0^{t\wedge
T_{n,\varepsilon}}-\ln\bigl(G_{nZ_n(r)}(\lambda)\bigr) a_n Z_n(r)^{\gamma}\,\d r
\]
converges to $\psi(\lambda)(t \wedge T_{\varepsilon})$ by
dominated convergence. Last, note that
\begin{eqnarray*}
\int_{\lfloor a_n\tau_n^{-1}(t\wedge T_{n,\varepsilon}) \rfloor
}^{a_n\tau_n^{-1}(t \wedge T_{n,\varepsilon})} -\ln\bigl( G_{X_n(\lfloor r
\rfloor)}(\lambda) \bigr) \,\d r &\leq& -\ln\bigl( G_{X_n (\lfloor a_n
\tau_n^{-1}(t) \rfloor)}
(\lambda) \bigr)\mathbf1_{\{ t<T_{n,\varepsilon}\}} \\
&=& -\ln\bigl( G_{nZ_n(t)} (\lambda) \bigr) \mathbf1_{\{ t<T_{n,\varepsilon}\}},
\end{eqnarray*}
since $X_n(\lfloor r \rfloor)$ is constant on the integration interval
and $a_n\tau_n^{-1}(t \wedge T_{n,\varepsilon})$ is an integer when $t
\geq T_{n,\varepsilon}$. The right-hand side in the inequality above
converges to $0$ as $n \rightarrow\infty$ since $nZ_n(t)>n
\varepsilon$ when $t< T_{n,\varepsilon}$ and $G_n(\lambda)
\rightarrow
1$ as $n \rightarrow\infty$. Finally, we have proved the convergence
(\ref{CVU}), hence~(a).
%This proves (a).

%Last, note that
% T_{n,\varepsilon})a_n]}^{\tau_n^{-1}(t_n \wedge
% T_{n,\varepsilon})a_n} -\ln( G_{X_n([r])}(\lambda) )
% ([\tau_n^{-1}(t_n)a_n] )}(\lambda) )
%G_{nZ_n(t_n)}(\lambda) ) \mathbf1_{\{ t_n<T_{n,\varepsilon}\}}
%(with the convention $0 \cdot\infty=0$), since $X_n([r])$ is
%constant on the integration interval and $\tau_n^{-1}(t_n\wedge
%T_{n,\varepsilon})a_n$ is an integer when $t_n \geq
%T_{n,\varepsilon}$. The right-hand side in the inequality above
%converges to $0$ as $n \rightarrow\infty$ since $nZ_n(t_n))>n
%convergence (\ref{CVU}), hence a).

The assertion (b) follows as a simple consequence of (a). By $(\ref
{MajorationMartingale})$ we have that, for each $t \geq0$, $(
M_n^{(\lambda)}(t \wedge T_{n,\varepsilon} ), n \geq n_0(\lambda)/
\varepsilon) $ is uniformly integrable. Together with the convergence
of (a), this is sufficient to deduce that the limit process $(Z'(t
\wedge
T_{\varepsilon})^{\lambda}\exp(\psi(\lambda) (t \wedge
T_{\varepsilon} )), t \geq0)$ is a martingale with respect to its
natural filtration. See \cite{EK}, Example 7, page 362.

We finally prove (c). Note that
\[
\bigl(Z'(t \wedge T_{\varepsilon})^{\lambda}\exp\bigl(\psi(\lambda) (t
\wedge
T_{\varepsilon} )\bigr), t \geq0\bigr) \mathop{\longrightarrow}\limits_{\varepsilon\rightarrow0
%, \varepsilon\notin\mathcal C_1 \cup\mathcal C_2
} \bigl( Z'(t)^{\lambda}\exp(\psi(\lambda) t ), t \geq0\bigr)
\]
for the Skorokhod topology. Besides, for each $t \geq0$ and
$\eps>0$, we have
\[
 Z'(t \wedge
T_{\varepsilon})^{\lambda}\exp(\psi(\lambda) (t \wedge
T_{\varepsilon} )) \leq\exp(\psi(\lambda) t ).
\]
As before, we can use
an argument of uniform integrability to conclude that\break
$(Z'(t)^{\lambda}\exp(\psi(\lambda) t ), t \geq0)$ is a martingale.
\end{pf*}

%is a power function, because the truncation argument is then useless
%to prove that $(Z^{\lambda}(t)\exp(\psi(\lambda) t ), t \geq0)$ is a
%martingale.

%%%%%%%%%%%%%%%%%%
%s4.3 ###
\subsection{Absorption times}
\label{Absorption}
%%%%%%%%%%%%%%%%%%

Recall that $A_n$ denotes the first time at which $X_n$ reaches the set
of absorbing states $\mathcal A$. To start with, we point out that
there is no loss of generality in assuming that $\mathcal A=\{0\}$.
Indeed, let $a_{\max}$ be the largest element of $\mathcal A$. If
$a_{\max} \geq1$, one can build a Markov chain $\tilde X_n$ starting
from $n$ and with transition probabilities $\tilde p_{i,j}=p_{i,j}$ for
$i \notin\mathcal A$ and all $j \geq0$, $\tilde p_{i,0}=1$ for $i \in
\mathcal A$, so that
\[
\tilde X_n(k)=X_n(k)\qquad \mbox{for } k \leq A_n \quad\mbox{and}\quad \tilde
X_n(k)=0\qquad \mbox{if } k > A_n.
\]%\vadjust{\goodbreak}
Clearly, this modified chain has a unique absorbing state, which is 0,
%If $\tilde{G}_n(\lambda)$ denotes the generating function
%associated with $(\tilde{p}_{n,k})$, as in (\ref{fonge}), then the
%functions $\tilde{G}_n$ and $\tilde{G}_n$ are equal for every
%$n>a_{\max}$. Therefore,
and the transition probabilities $(\tilde{p}_{n,k})$ satisfy (H)
if and only if $(p_{n,k})$ do. Besides, the first time $\tilde A_n$ at
which $\tilde X_n$ reaches~0 is clearly either equal to $A_n$ or to
$A_n+1$.
%Theorem
%a_{\max})$, it will also hold for $(A_n, n \geq a_{\max})$.
Moreover, constructing $\tilde Y_n$ from $\tilde X_n$ as $Y_n$ is
defined from $X_n$ in Section \ref{secscal}, we see that $\sup_{t
\geq
0} | \tilde Y_n(t)-\tilde Y_n(t) | \leq a_{\max}/n$. This is enough to
see that the convergence in distribution as $n \rightarrow\infty$ of
$(\tilde A_n/a_n,\tilde Y_n)$ entails that of $(A_n/a_n,Y_n)$ towards
the same limit. This in turn entails the convergence in distribution of
$(A_n/a_n,Y_n,Z_n)$ to the required limit, using a part of the proof of
\cite{EK}, Theorem 1.5, Chapter 6, as already mentioned at the
beginning of Section \ref{Limit}. In conclusion, if the convergence of
Theorem \ref{cvdeathtimes} is proved for the sequence $(\tilde A_n/a_n,
\tilde Y_n, \tilde Z_n)$, $n \geq0$, it will also hold for $(A_n/a_n,
Y_n, Z_n)$, $n \geq0$, with the same distribution limit. In the
following, we will therefore additionally suppose that $a_{\max}=0$,
that is,
%
%e28 ###
\begin{equation} \label{Hopla}
\mathcal A=\{0\} ,
\quad\mbox{or equivalently,}\quad p_{n,n}<1 \qquad\mbox{for every }n\geq1.
\end{equation}
We now set out a preliminary lemma that we will use for the proof of
Theorem \ref{cvdeathtimes}.

\begin{lemma}
\label{majorationZn} For every $\lambda>0$ and $\beta>0$, there exists
some finite constant $c_{\lambda,\beta}>0$ such that for all $n \in
\Z_+$ and all $t \geq0$,
%
%e29 ###
\begin{equation}\label{eq:8}
Z_n(t)^{\lambda} \leq\frac{c_{\lambda,\beta} M_n^{(\lambda)}(t) +1
}{t^{\beta}} ,
\end{equation}
where the processes $M_n^{(\lambda)}$ are the martingales defined in
(\ref{martingaleMn}). Consequently,
\[
\mathbb E[Z_n(t)^{\lambda}] \leq\frac{c_{\lambda,\beta} +1
}{t^{\beta}}.
\]
\end{lemma}

In the cases where $n^{-\gamma} a_n \rightarrow\ell\in(0,\infty)$,
our proof can be adapted to get the following stronger result: There
exists some finite constant $c_{\lambda}$ such that for all $n \in
\mathbb Z_+$ and all $t \geq0$, $ Z_n(t)^{\lambda} \leq
M_n^{(\lambda)}(t) \exp(c_{\lambda}(1-t)), $ and consequently, $
\E[Z_n(t)^{\lambda}] \leq\exp(c_{\lambda}(1-t)). $

\begin{pf}
Fix $\lambda>0$ and $\beta>0$. For a given
$n,t$, if
$Z_n(t)^{\lambda}\leq t^{-\beta}$, then obviously (\ref{eq:8}) is satisfied,
irrespective of any choice of $c_{\lambda,\beta}$. So we assume
that $Z_n(t)^{\lambda}>t^{-\beta}$, and in particular, $Z_n(t)>0$. By
(\ref{martingaleMn}), we have
\[
Z_n(t)^{\lambda} = M^{(\lambda)}_n(t) \exp\biggl(\int_0^{\lfloor
a_n\tau_n^{-1}(t)\rfloor} \ln\bigl( G_{X_n(\lfloor r\rfloor)}(\lambda)
\bigr)\,\d{r} \biggr).
\]
Note that $X_n(\lfloor r \rfloor) \geq1$ as soon as $r \leq\lfloor
a_n\tau_n^{-1}(t)\rfloor$ and $Z_n(t) > 0$. Moreover, under the
assumption~(\ref{Hopla}), we have $ \ln( G_{n}(\lambda) )
\leq-c_4(\lambda)/a_n <0$ for every $n \geq1$ by (\ref{ceq3}). Hence,
for all $\varepsilon>0$,
\begin{eqnarray*}
&&\int_0^{\lfloor a_n\tau_n^{-1}(t)\rfloor} \ln\bigl( G_{X_n(\lfloor r\rfloor)}(\lambda) \bigr) \,\mathrm{d}r
\\
&&\quad\leq
\int_0^{\lfloor a_n\tau_n^{-1}(t) \rfloor} \frac{-c_4(\lambda)}{a_{X_n(\lfloor r\rfloor)}} \,\mathrm d r
=
\int_0^{a_n\tau_n^{-1}(t)}\frac{-c_4(\lambda)}{a_{X_n(\lfloor r\rfloor)}} \,\mathrm d r-\int_{\lfloor a_n\tau_n^{-1}(t)\rfloor}^{a_n\tau_n^{-1}(t) }
\frac{-c_4(\lambda)}{a_{X_n(\lfloor r\rfloor)}} \,\mathrm d r %\vadjust{\goodbreak}
\\
&&\quad\mathop{\leq}\limits_{\mathrm{by}\ (\ref{chvar})}
- c_4(\lambda) \int_0^t\frac{a_n}{a_{nZ_n(r)}} Z_n(r)^{\gamma} \,\mathrm dr+\frac{c_4(\lambda)}{\inf_{k \geq0}a_k}
\\
&&\quad\mathop{\leq}\limits_{\mathrm{by}\ (\ref{ceq2})}
- c_4(\lambda) c'_1(\varepsilon) \int_0^t Z_n(r)^{\varepsilon} \,\mathrm dr+\frac{c_4(\lambda)}{\inf_{k\geq0}a_k}
\leq
- c_4(\lambda) c'_1(\varepsilon) t Z_n(t)^{\varepsilon} + \frac{c_4(\lambda)}{\inf_{k \geq0 }a_k}.
\end{eqnarray*}
Since $Z_n(t)^{\lambda}>t^{-\beta}$, we have, taking
$\varepsilon=\lambda/2\beta$, the existence of a finite constant
$c_{\lambda,\beta}$, independent of $n$ and $t$, such that
\[
Z_n(t)^{\lambda}\leq M_n^{(\lambda)}(t) \exp\Bigl(-
c_4(\lambda)c'_1(\lambda/2\beta) t^{1/2} +c_4(\lambda)/\inf_{k\geq
0}a_k\Bigr) \leq c_{\lambda,\beta}M_n^{(\lambda)}(t)/t^{\beta} ,
\]
giving the result.
\end{pf}

\begin{pf*}{Proof of Theorem \ref{cvdeathtimes}}
Notice that the
first time at which $Y_n$ reaches 0 is
\[
\int_0^{\infty} Z_n(r)^{\gamma}\,\mathrm dr=\sigma_n=A_n/a_n
\]
using (\ref{sigma}) for the first equality and (\ref{Hopla}) for
the second equality. The previous lemma ensures that $\sup_{n \geq1}
\mathbb E[\sigma_n] <\infty$, which implies that the sequence
$(\sigma_n, n \geq1 )$ is tight. In turn, this implies that the
sequence $((Y_n,Z_n,\sigma_n), n \geq1 )$ is tight.

The proof of Theorem \ref{cvdeathtimes} will therefore be completed if
we prove the uniqueness of possible limiting distributions of
$((Y_n,Z_n,\sigma_n), n \geq1 )$ along a subsequence. In that aim,
consider a strictly increasing sequence of integers $(n_k,k\geq0)$ such
that the sequence $((Y_n,Z_n,\sigma_n), n \geq1)$ converges in
distribution along $(n_k)$ to a limit $(Y',Z',\sigma')$. By Proposition~\ref{sec:scaling-limits-non-2}, $(Y',Z')$ has same distribution as
$(Y,Z)$, so by abuse of notations, for simplicity, we write $(Y,Z)$
instead of $(Y',Z')$. Our goal is to show that $\sigma'$ is the
extinction time $\sigma=\sigma_Y=\int_0^{\infty} Z(r)^{\gamma}
\,\mathrm
dr$, with the notations of Section \ref{tc}.

By the Skorokhod representation theorem, we may suppose that the
convergence of $(Y_n,Z_n,\sigma_n)$ to $(Y,Z,\sigma')$ is almost sure.
It is then immediately checked that a.s.,
\[
\sigma'=\liminf_{n\to\infty}\inf\{t\geq0 \dvtx Y_n(t)=0\}\geq\inf\{
t\geq
0\dvtx Y(t)=0\}=\sigma,
\]
so in order to show that $\sigma=\sigma'$ a.s., it suffices to check
that $\E[\sigma']\leq\E[\sigma]$. To see this, note that the
convergence in the Skorokhod sense implies that a.s., for a.e. $t$,
$Z_n(t) \rightarrow Z(t)$ and therefore, by Fubini's theorem, that for
a.e. $t$, $Z_n(t) \rightarrow Z(t)$ a.s. We then obtain that for a.e.
$t$, $Z_n(t)^{\gamma} \rightarrow Z(t)^{\gamma}$ a.s., and since all
these quantities are bounded by $1$, we have, by dominated convergence,
that for a.e. $t$, $\mathbb E[Z_n(t)^{\gamma}] \rightarrow\mathbb
E[Z(t)^{\gamma}]$. Then, again by dominated convergence, using Lemma
\ref{majorationZn}, we get $\int_0^{\infty} \mathbb E[Z_n(r)^{\gamma}]
\,\mathrm dr \rightarrow\int_0^{\infty} \mathbb E[Z(r)^{\gamma}]
\,\mathrm
dr<\infty$. Hence, by Fubini's theorem
\[
\mathbb E \biggl[ \int_0^{\infty} Z_n(r)^{\gamma} \,\mathrm dr \biggr] \mathop{\longrightarrow}\limits_{n
\rightarrow\infty}\mathbb E \biggl[ \int_0^{\infty}
Z(r)^{\gamma} \,\mathrm d r \biggr].
\]
But,
%again
by Fatou's lemma,
\[
\E[\sigma'] \leq\liminf_n \E\biggl[ \int_0^{\infty} Z_n(r)^{\gamma}
\,\mathrm
dr \biggr] = \mathbb E \biggl[ \int_0^{\infty} Z(r)^{\gamma} \,\mathrm d r
\biggr]=\E[\sigma] ,
\]
%
% Together with (\ref{Fatou1}), and since these expectations are
%finite, this implies that $\sigma'=\int_0^{\infty} Z(r)^{\gamma}
as wanted. This shows that $(Y_n,Z_n,\sigma_n)$ converges in
distribution (without having to take a subsequence) to $(Y,Z,\sigma)$,
which gives the first statement of Theorem \ref{cvdeathtimes}.

%We can then conclude with Theorem \ref{CvChain} that
%the distribution of $(Y,Z,\sigma)$ is then necessarily that of
%$(\exp(-\xi_\tau),\exp(-\xi),\int_0^{\infty}\exp(-\gamma\xi_r)\mathrm
%dr )$, where $\xi$ is a subordinator with Laplace exponent $\psi$ and
%$\tau$ the time change defined in the statement of Theorem

%With a slight abuse of notations, we set $\sigma=\int_0^{\infty}
Since $(Y_n,Z_n,\sigma_n)$ converges in distribution to $(Y,Z,\sigma)$,
by using Skorokhod's representation theorem, we assume that the
convergence is almost-sure. Note that the above proof actually entails
the convergence of moments of order 1, $\mathbb E[\sigma_n]
\rightarrow
\mathbb E[\sigma]$. We now want to prove the convergence of moments of
orders $u \geq0$. It is well known (see \cite{CarmonaPetitYor},
Proposition~3.3) that the random variable $\sigma$ has positive moments
of all orders and that its moment of order $p \in\N$ is equal to
$p!/\prod_{j=1}^p \psi(\gamma j)$. Let $u \geq0$. Since $\sigma_n
\rightarrow\sigma$ a.s., if we show that $\sup_{n
\geq1}\mathbb E[\sigma_n^{p}]<\infty$ for some $p>u,$ then
$((\sigma_n-\sigma)^u,n\geq0)$ will be uniformly integrable, entailing
the convergence of $\E[|\sigma_n-\sigma|^u]$ to 0. So fix $p>1$,
consider $q$ such that $p^{-1}+q^{-1}=1$ and use H\"{o}lder's
inequality to get
\[
\int_1^{\infty} Z_n(r)^{\gamma} \,\mathrm dr \leq\biggl( \int_1^{\infty}
Z_n(r)^{\gamma p} r^{2p/q} \,\mathrm dr \biggr)^{1/p} \biggl( \int_1^{\infty} r^{-2}
\,\mathrm dr\biggr)^{1/q}.
\]
Together with Lemma \ref{majorationZn} this implies that
\[
\sup_{n \geq1} \mathbb E \biggl[\biggl( \int_1^{\infty} Z_n(r)^{\gamma}
\,\mathrm dr
\biggr)^p \biggr]<\infty,
\]
which, clearly, leads to the
required $\sup_{n \geq1}\mathbb E[\sigma_n^{p}]<\infty$.
\end{pf*}

%s4.4 ###
\subsection{\texorpdfstring{Proof of Proposition
\protect\ref{sec:scaling-limits-non-3}}{Proof of Proposition 1}}

Consider a probability measure $\mu$ on $[0,1]$, a real number
$\gamma>0$ and a function $\ell\dvtx \mathbb R_+ \rightarrow(0,\infty)$
slowly varying at $\infty$. Then set $a_n= n^{\gamma}\ell(n)$, let
$\gamma'$ be such that $\max(1, \gamma) < \gamma' < \gamma+1$ and
assume that $n$ is large enough so that $n^{\gamma'-1}<a_n \leq
n^{\gamma'}$. For such an $n$ and $0\leq k\leq n-1$, set
\begin{eqnarray*}
p_{n,k} &=& a_n^{-1} \int_{[0,1-a_n^{-1})} \pmatrix{n\cr k} x^k(1-x)^{n-k-1}
\mu(\mathrm dx)
\\
&&{}+ n^{1-\gamma'} \mu(\{1\}) \mathbf1_{\{k= n-\lfloor
n^{\gamma'}/a_n\rfloor\}},\qquad 0\leq k \leq n-1.
\end{eqnarray*}
Clearly, these quantities are non-negative and
\[
\sum_{k=0}^{n-1} p_{n,k}=a_n^{-1}
\int_{[0,1-a_n^{-1})}(1-x^n)(1-x)^{-1} \mu(\mathrm dx)
+n^{1-\gamma'}\mu(\{1\}) \leq\mu([0,1)) + \mu(\{ 1\})=1.
\]
Let $p_{n,n} = 1-\sum_{k=0}^{n-1} p_{n,k},$ in order to define a
probability vector $(p_{n,k},0 \leq k \leq n)$ on $\{0,1,\ldots,n \}$.
Now, for any continuous test function $f \dvtx [0,1] \rightarrow\mathbb
R_+$,
\begin{eqnarray*}
a_n \sum_{k=0}^{n} f\biggl(\frac{k}{n}\biggr) \biggl(1-\frac{k}{n}\biggr) p_{n,k}&=&
\int_{[0,1-a_n^{-1})}\sum_{k=0}^{n-1} f\biggl(\frac{k}{n}\biggr) \biggl(1-\frac{k}{n}\biggr)
\pmatrix{n\cr k}
x^k(1-x)^{n-k-1}\mu(\mathrm dx) \\
&&{}+ a_n n^{-\gamma'} \lfloor n^{\gamma'}/a_n\rfloor f(1-n^{-1}\lfloor
n^{\gamma'}/a_n\rfloor) \mu(\{1\}).
\end{eqnarray*}
The term involving $\mu(\{1\})$ clearly converges to $f(1)\mu(\{1\})$
since $f$ is continuous. For the other term, note that for $B_{n,x}$ a
binomial random variable with parameters $n,x$,
\[
\sum_{k=0}^{n-1} f\biggl(\frac{k}{n}\biggr) \biggl(1-\frac{k}{n}\biggr) \pmatrix{n\cr k}
x^k(1-x)^{n-k}=\mathbb E \biggl[ f\biggl(\frac{B_{n,x}}{n}\biggr)\biggl(1-\frac{B_{n,x}}{n} \biggr)\biggr],
\]
which converges to $f(x)(1-x)$ as $n \rightarrow\infty$ and is
bounded on $[0,1]$ by a constant times $(1-x)$ since $f$ is bounded.
Hence by dominated convergence,
\[
a_n \sum_{k=0}^{n} f\biggl(\frac{k}{n}\biggr) \biggl(1-\frac{k}{n}\biggr) p_{n,k}
\rightarrow\int_{[0,1]}f(x) \mu(\mathrm dx).
\]

%%%%%%%%%%%%%%%%%%%%%%%%%%
%%%%%%%%%%%%%%%%%%%%%%%%%%

%%%%%%%%%%%%%%%%%%%%%%%%%%%%%%%%
%s5 ###
\section{Scaling limits of random walks with a barrier}
\label{RW}
%%%%%%%%%%%%%%%%%%%%%%%%%%%%%%%%

Recall that %$(\zeta_i,i \geq1)$ is an sequence of
%i.i.d. $\Z_+$-valued random variables with distribution
$(q_k,k\geq0)$ is a probability distribution satisfying $q_0<1$, as
well as the definition of the random walk with a barrier model $X_n = n
-S^{(n)}$ and notation from Section \ref{sec:random-walk-with}. In the
following, $n$ will always be implicitly assumed to be large enough so
that $\overline{q}_n<1$.\looseness=-1

%Of course, $S_k \rightarrow\infty$ a.s., as $k \rightarrow
%%\infty$. Our goal here is to consider modified versions of this
%random walk, by forcing it not to exceed the barrier $n$, in the
%following different manners.

%In the following, we set $i_0:=\inf\{k \geq0: q_k>0 \}$. Note then
%that $n-U_n$, $n-V_n$ and $n-W_n$ are non-increasing Markov chains
%starting from $n$ and with transition
%probabilities \begin{enumerate} \item[$\bullet$]
%$p^{(U)}_{i,j}=q_{i-j}$ if $1 \leq j \leq i$; $p^{(U)}_{i,j}=\sum_{k
%and $p_{n,0}=0$ \item[$\bullet$] if $i>i_0$,
%$p^{(W)}_{i,j}=q_{i-j}/(q_{i_0}+...+q_{i-1})$ for $1 \leq j \leq
%i-i_0$, 0 otherwise; and if $i \leq i_0$,
%$p_{i,i}=1$. \end{enumerate}

\begin{pf*}{Proof of Theorem \ref{ThRW}}
Let us first prove (i).
We assume that $\overline{q}_n = n^{-\gamma}\ell(n)$, where $(\ell
(x),x\geq
0)$ is slowly varying at $\infty$. We want to show that (H) is
satisfied, with $a_n=1/\overline{q}_n$ and
$\mu(\d{x})=\gamma(1-x)^{-\gamma}\,\d{x}\,\ind_{\{0<x<1\}}$. From
this, the
conclusion follows immediately.

Using the particular form of the transition probabilities (\ref{eq:5}),
it is sufficient to show that for every function $f$ that is
continuously differentiable on $[0,1]$,
%
%e30 ###
\begin{equation}\label{eq:10}
\frac{1}{\overline{q}_n}\sum_{k=0}^n
\frac{q_{n-k}}{\sum_{i=0}^nq_i}\biggl(1-\frac{k}{n}\biggr)f\biggl(\frac
{k}{n}\biggr) \longrightarrow\gamma\int_0^1f(x)(1-x)^{-\gamma}\,\d{x} .
\end{equation}
Let $g(x)=xf(1-x)$. By Taylor's expansion, we have, for every $x\in
(0,1)$, $g((x+\frac{1}{n})\wedge1) - g(x)=g'(x)/n+\eps_n(x)/n$, where
$\sup_{x\in[0,1]}\eps_n(x)$ converges to $0$ as $n\to\infty$.
Therefore, since $g(0)=0$,
\begin{eqnarray*}
\frac{1}{\overline{q}_n}\sum_{k=0}^n\frac{q_{n-k}}{\sum_{i=0}^nq_i}\biggl(1-\frac{k}{n}\biggr)f\biggl(\frac{k}{n}\biggr)
&=&
\frac{1}{\overline{q}_n(1-\overline{q}_n)}\sum_{k=0}^nq_k g\biggl(\frac{k}{n}\biggr)
\\%[-3pt]
&=&
\frac{1}{\overline{q}_n(1-\overline{q}_n)}\sum_{k=0}^{n-1}\overline{q}_k\biggl(g\biggl(\frac{k+1}{n}\biggr)-g\biggl(\frac{k}{n}\biggr)\biggr)-\frac{g(1)}{1-\overline{q}_n}
\\%[-3pt]
&=&
\frac{1}{n\overline{q}_n(1-\overline{q}_n)}\sum_{k=0}^{n-1}\overline{q}_k\biggl( g'\biggl(\frac{k}{n}\biggr)+\eps_n\biggl(\frac{k}{n}\biggr)\biggr)-\frac{g(1)}{1-\overline{q}_n} .
\end{eqnarray*}
Because of the uniform convergence of $\eps_n$ to $0$, this is
equivalent as $n\to\infty$ to
\[
\frac{1}{n(1-\overline{q}_n)}\sum_{k=0}^{n-1} \frac{\overline
{q}_k}{\overline{q}_n}
g'\biggl(\frac{k}{n}\biggr)-\frac{g(1)}{1-\overline{q}_n}
\mathop{\longrightarrow}\limits_{n\to\infty}\int_0^1x^{-\gamma}g'(x)\,\d
{x}-g(1)
\]
by a simple use of the uniform convergence theorem for regularly
varying functions (\cite{BGT}, Theorem 1.5.2). Integrating by parts,
the latter integral is the right-hand side of (\ref{eq:10}).

Statement (ii) is even simpler.
% We use Proposition \ref{sec:preliminaries}, and
Fix $\lambda>0$. For all $k\geq1$, it holds that $n(1-(
1-k/n)^\lambda)$ converges to $\lambda k$ as $n \rightarrow\infty$.
Moreover, $n (1- ( 1-k/n )^{\lambda} ) \leq\max(1,\lambda) k$. Hence,
when $m=\sum_{k=0}^{\infty}kq_k<\infty$, we have by dominated
convergence that
%$$ n \sum_{k=0}^{n}q_k (1- ( 1-\frac{k}{n}
%)^{\lambda} )\rightarrow\lambda m.
%$$ %Besides, $\sum_{k \geq n} q_k \leq n^{-1} \sum_{k \geq n}
%kq_k=\circ(n^{-1})$.
%t is thus clear that
%
\[
n\Biggl(1-\sum_{k=0}^n p_{n,k} \biggl( \frac{k}{n} \biggr)^{\lambda
}\Biggr)
\mathop{\longrightarrow}\limits_{n\to\infty} \lambda m .
\]
We conclude by Theorems \ref{CvChain} and \ref{cvdeathtimes} and
Proposition~\ref{sec:preliminaries}.
\end{pf*}

Let us now consider some variants of the random walk with a barrier.
The results below recover and generalize results of \cite{IM08}. Let
$(\zeta_i,i\geq1)$ be an i.i.d. sequence with distribution
$(q_n,n\geq0)$. Set $S_0=0$ and
\[
S_k=\sum_{i=1}^k\zeta_i ,\qquad k\geq1 ,
\]
for the random walk
associated with $(\zeta_i,i\geq1)$. We let $\widetilde{S}^{(n)}_k = n
\wedge
S_k,k\geq0$ be the walk truncated at level $n$. We also define
$\widehat{S}_k^{\,(n)}$ recursively as follows: $\widehat
{S}^{\,(n)}_0=0$, and given
$\widehat{S}^{\,(n)}_k$ has been defined, we let
\[
\widehat{S}^{\,(n)}_{k+1}=\widehat{S}^{\,(n)}_k+\zeta_{k+1}\ind_{\{
\widehat
{S}^{(n)}_k+\zeta_{k+1}\leq n\}} .
\]
In other words, the process $\widehat{S}^{\,(n)}$ evolves as
$S$, but ignores the jumps that would bring it to a level higher than
$n$. This is what is called the random walk with a barrier in
\cite{IM08}. However, in the latter reference, the authors assume that
$q_0=0$ and therefore really consider the variable $A_n$ associated
with $X_n$ as defined above, as they are interested in the number of
\textit{strictly positive} jumps that $\widehat{S}^{\,(n)}$
accomplishes before
attaining its absorbing state. See the forthcoming Lemma
\ref{sec:scal-limits-rand-1} for a proof of the identity in
distribution between $A_n$ and the number of strictly positive jumps
of $\widehat{S}^{\,(n)}$ when $q_0=0$.

The processes $\widetilde{X}_n = n -\widetilde{S}^{(n)}$ and
$\widehat{X}_n = n
-\widehat{S}^{\,(n)}$ are non-increasing Markov chains with transition
probabilities given by
\[%\label{eq:12}
\widetilde{p}_{i,j}=q_{i-j}+\ind_{\{j=0\}}\overline{q}_i ,\qquad
\widehat{p}_{i,j}=q_{i-j}+\ind_{\{j=i\}}\overline{q}_i ,\qquad 0\leq j\leq
i .
\]
We let $\widetilde{A}_n$ and $\widehat{A}_n$ be the respective
absorption times. By
an argument similar to that in the above proof, it is easy to show that
when $\overline{q}_n$ is of the form $n^{-\gamma}\ell(n)$ for some
$\gamma
\in
(0,1)$ and slowly varying function $\ell$, then (H) is satisfied
for these two models, with sequence $a_n=1/\overline{q}_n$ and measures
\[
\widetilde{\mu}=\delta_0 + \mu=\delta_0+\gamma(1-x)^{-\gamma}\,\d{x}\,\ind_{\{0<x<1\}}
,\qquad\widehat{\mu}=\mu=\gamma(1-x)^{-\gamma}\,\d{x}\,\ind_{\{0<x<1\}} .
\]
Consequently, we obtain the joint convergence of
$\widehat{Y}_n=(\widehat{X}_n(\lfloor t/\overline{q}_n\rfloor
)/n,t\geq0)$ and
$\overline{q}_n\widehat{A}_n$ to the same distributional limit as
$(Y_n,\overline{q}_nA_n)$ as in (i), Theorem \ref{ThRW}, with the obvious
notation for $Y_n$. In the same way, $(\widetilde{Y}_n(\lfloor
t/\overline{q}_n\rfloor),t\geq0)$ and $\widetilde{A}_n$ converge to
the limits
involved in Theorems \ref{CvChain} and \ref{cvdeathtimes}, but this
time, using a killed subordinator $\xi^{(\mathtt{k})}$ with Laplace
exponent
\[
\psi^{(\mathtt{k})}(\lambda)=\psi(\lambda)+1=1+ \int_0^\infty
(1-\mathrm{e}^{-\lambda y}) \frac{\gamma \mathrm{e}^{-y}\,\d{y}}{(1-\mathrm{e}^{-y})^{\gamma+1}}
,\qquad\lambda\geq0 .
\]
If $\xi$ is a subordinator with Laplace
exponent $\psi$, and if $\be$ is an exponential random variable with
mean $1$, independent of $\xi$, then
$\xi^{(\mathtt{k})}(t)=\xi(t)+\infty\ind_{\{t\geq\be\}},t\geq0$
is a
killed subordinator with Laplace exponent $\psi^{(\mathtt{k})}$.

In fact, we have a joint convergence linking the processes
$X_n,\widetilde{X}_n,\widehat{X}_n$ together. Note that the three can
be joined
together in a very natural way, by building them with the same
variables $(\zeta_i,i\geq1)$. This is obvious for $\widetilde{X}_n$ and
$\widehat{X}_n$, by construction. Now, a process with the same distribution
as $X_n$ can be constructed simultaneously with $\widehat{X}_n$ by a simple
time change, as follows.

\begin{lemma}\label{sec:scal-limits-rand-1}
Let $T^{(n)}_0=0$, and recursively, let
\[
T^{(n)}_{k+1}=\inf\bigl\{i> T^{(n)}_k\dvtx \widehat{S}^{\,(n)}_{i-1}+\zeta
_{i}\leq
n\bigr\} .
\]
%
% Let $$R_n(k)=\sum_{i=1}^{k}\ind_{\{\wh{S}^{(n)}_{i-1}+\zeta_{i}\leq
%n\}} , R_n^{-1}(k)=\max\{i\geq0:R_n(i)= k\}$$
Then the process $(\widehat{X}_n(T^{(n)}_k),k\geq0)$ has same distribution
as $X_n$, with the convention that $\widehat{X}_n(\infty)=\lim_{k
\rightarrow\infty} \widehat{X}_n(k)$.
\end{lemma}

\begin{pf}
We\vspace*{-2pt} observe that the sequence
$(\zeta_{T^{(n)}_k},k\geq1)$ is constructed by rejecting elements
$\zeta_i$ such that $\widehat{S}^{\,(n)}_{i-1}+\zeta_i>n$, so by a simple
recursive argument, given $\widehat{X}_n(T^{(n)}_k)$, the random variable
$\zeta_{T^{(n)}_{k+1}}$ has the same distribution as a random variable
$\zeta$ with distribution $q$ conditioned on
$\widehat{X}_n(T^{(n)}_k)+\zeta\leq n$. This is exactly the
definition of
$S^{(n)}$.
\end{pf}

In the following statement, we assume that $X_n,\widetilde
{X}_n,\widehat{X}_n$ are
constructed jointly as above. We let $\xi$ be a subordinator with
Laplace exponent $\psi$ as in the statement of (i) in Theorem
\ref{ThRW}. Let $\xi^{(\mathtt{k})}$ be defined as above, using an
independent exponential variable $\be$. Let $\tau$ be the time change
defined as in the Theorem \ref{ThRW}, and let $\tau^{(\mathtt{k})}$ be
defined similarly from $\xi^{(\mathtt{k})}$. Let
\[
Y=\bigl(\exp\bigl(-\xi_{\tau(t)}\bigr),t\geq0\bigr) ,\qquad
\widetilde{Y}=\bigl(\exp\bigl(-\xi^{(\mathtt{k})}_{\tau^{(\mathtt
{k})}(t)}\bigr),t\geq
0\bigr) .
\]

\begin{proposition}\label{sec:scal-limits-rand}
Under the same hypotheses as in {\rm(i)}, Theorem \ref{ThRW}, the
following convergence in distribution holds in $\mathcal{D}^3$:
\[
(Y_n,\widetilde{Y}_n,\widehat{Y}_n)\mathop{\stackrel{(d)}{\longrightarrow}}_{n\to\infty}(Y,\widetilde{Y},Y) ,
\]
and jointly,
\[
\overline{q}_n(A_n,\widetilde{A}_n,\widehat{A}_n)\mathop{\stackrel{(d)}{\longrightarrow}}_{n\to\infty}
\biggl(\int_0^{\infty}\mathrm{e}^{-\gamma\xi_t}\,\d
{t},\int_0^{\be}\mathrm{e}^{-\gamma\xi_t}\,\d{t},\int_0^{\infty}\mathrm{e}^{-\gamma
\xi_t}\,\d{t}\biggr).
\]
\end{proposition}

\begin{pf*}{Proof (sketch)} The convergence of one-dimensional
marginals holds by the above discussion. Let
$(Y^{(1)},Y^{(2)},Y^{(3)},\sigma^{(1)},\sigma^{(2)},\sigma^{(3)})$
be a
limit in distribution of the properly rescaled $6$-tuple
$(X_n,\widetilde{X}_n,\widehat{X}_n,A_n,\widetilde{A}_n,\widehat
{A_n})$ along some subsequence.
These variables are constructed by three subordinators,
$\xi^{(1)},\xi^{(2)},\xi^{(3)}$, with the\vspace*{1pt} same law as
$\xi,\xi^{(\mathtt{k})},\xi,$ respectively. Now, we use the obvious
fact that $\widetilde{X}_n\leq X_n\leq\widehat{X}_n$. Taking limits,
we have
$Y^{(2)}\leq Y^{(1)}\leq Y^{(3)}$ a.s. Taking expectations, using the
fact that $Y^{(1)}$ and $Y^{(3)}$ have the same distribution and using
the fact that these processes are c\`{a}dl\`{a}g, we obtain that
$Y^{(1)}=Y^{(3)}$ a.s. Similarly, $\sigma^{(1)}=\sigma^{(3)}\geq
\sigma^{(2)}$ a.s., and $\sigma^{(2)}$ is the first time where
$Y^{(2)}$ attains $0$ (which\vspace*{1pt} is done by accomplishing a negative jump).
Moreover, we have $\widetilde{X}_n(k)= X_n(k)= \widehat{X}_n(k)$ for
every $k<
\widetilde{A}_n$. By passing to the limit, we obtain that
$Y^{(1)}=Y^{(2)}=Y^{(3)}$ a.s. on the interval $[0,\sigma^{(2)}]$. This
shows that $\xi^{(1)}=\xi^{(3)}$ and that
$\xi^{(1)}=\xi^{(2)}=\xi^{(3)}$ on the interval where $\xi^{(2)}$ is
finite. Since $\xi^{(2)}$ is a killed subordinator, this completely
characterizes the distribution of $(\xi^{(1)},\xi^{(2)},\xi^{(3)})$ as
that of $(\xi,\xi^{(\mathtt{k})},\xi)$, and this allows us to conclude.
Details are left to the reader.
\end{pf*}

%%%%%%%%%%%%%%%%%%%%%%%%%%%%%%%%
%s6 ###
\section{Collisions in $\Lambda$-coalescents}
\label{Coalescents}
%%%%%%%%%%%%%%%%%%%%%%%%%%%%%%%%
We now prove Theorem \ref{theolambda}. Using Theorems \ref{CvChain} and
\ref{cvdeathtimes}, all we have to check is that the hypothesis~(H) is satisfied with the parameters $a_n=\int_{[1/n,1]}
x^{-2}\Lambda(\mathrm d x)$, $n \geq1$, and $\psi$\vspace*{1pt} defined by
(\ref{eq:7}). This is an easy consequence of the following Lemmas
\ref{lemlamb1} and \ref{lemlamb2}. We recall that the transition
probabilities of the Markov chain $(X_n(k), k \geq0)$, where $X_n(k)$
is the number of blocks after $k$ coalescing events when starting with
$n$ blocks, are given by (\ref{defprobacoal}).
%$$
%p_{n,k}=\frac{g_{n,k}}{g_n}=\frac{1}{g_n}\binom{n}{k-1} \int_{[0,1]}
%x^{n-k-1}(1-x)^{k-1}
%$$
%where $g_n$ is the normalizing constant.

\begin{lemma}
\label{lemlamb1}
Assume that $\Lambda(\{ 0\})=0 $ and that $u
\rightarrow\int_{[u,1]} x^{-2}\Lambda(\mathrm d x)$ varies regularly
at 0 with index $-\gamma$, $\gamma\in(0,1)$. Then,
\[
g_n \sim\Gamma(2-\gamma) \int_{[1/n,1]} x^{-2}\Lambda(\mathrm d x)
\qquad\mbox{as } n \rightarrow\infty.
\]
\end{lemma}

\begin{pf} First note that
\begin{eqnarray*}
g_n&=& \int_{(0,1]} \bigl(1-(1-x)^n-n(1-x)^{n-1}x \bigr)x^{-2} \Lambda(\mathrm
dx) \\ &=& I_n-J_n,
%&=& \int_{(0,\infty]} (1-\exp(-nx)-n\exp(-x(n-1))(1-\exp(-x))
%)(1-\exp(-x))^{-2}\tilde\Lambda(\mathrm dx)
\end{eqnarray*}
where, defining by $\tilde\Lambda$ the push-forward of $\Lambda$ by
the mapping $x \mapsto-\log(1-x)$,
\begin{eqnarray*}
I_n&=& \int_{(0,\infty]} \bigl(1-\exp(-nx)-n\exp(-xn)x\bigr)\bigl(1-\exp
(-x)\bigr)^{-2}\tilde
\Lambda(\mathrm dx),\\[-3pt]
J_n&=& \int_{(0,\infty]} \bigl(n\exp\bigl(-x(n-1)\bigr)\bigl(1-\exp(-x)-x\exp(-x)\bigr)
\bigr)\bigl(1-\exp(-x)\bigr)^{-2}\tilde\Lambda(\mathrm dx).
\end{eqnarray*}
The integrand in the integral $J_n$ converges to 0 as $n \rightarrow
\infty$, for all $x \in(0,\infty]$. And, clearly, there exists some
finite constant $C$ such that for all $ x \in(0,\infty]$, and all $n
\geq1$,
\begin{eqnarray*}
&&\bigl| n\exp\bigl(-x(n-1)\bigr)\bigl(1-\exp(-x)-x\exp(-x)\bigr) \bigl(1-\exp(-x)\bigr)^{-2} \bigr|
\\[-3pt]
&&\quad
\leq C \bigl(1-\exp(-x)\bigr)^{-1}.
\end{eqnarray*}
Hence, by dominated convergence, $ J_n \rightarrow0 \mbox{ as } n
\rightarrow\infty. $ Next, $I_n$ can be rewritten as
\begin{eqnarray*}
I_n&=&\int_{(0,\infty]} \biggl(\int_{(0,x]} n^2u\exp(-nu)\,\mathrm du\biggr)
\bigl(1-\exp(-x)\bigr)^{-2}\tilde\Lambda(\mathrm dx) \\[-3pt]
&=&
n^{2}\int_{(0,\infty)}\exp(-nu)u
\biggl(\int_{[1-\exp(-u),1]}x^{-2}\Lambda(\mathrm dx)\biggr) \,\mathrm du.
\end{eqnarray*}
Since $ \int_{[u,1]} x^{-2}\Lambda(\mathrm d x)$ varies regularly as
$u\rightarrow0$ with index $-\gamma$,
\[
u \int_{[1-\exp(-u),1]}x^{-2}\Lambda(\mathrm dx)\mathop{\sim}\limits_{u\rightarrow0} u \int_{[u,1]}x^{-2}\Lambda(\mathrm dx)
\]
and these functions vary regularly at 0 with index $1-\gamma$. It is
then standard that
\[
\int_{[0,t]}u \biggl( \int_{[1-\exp(-u),1]}x^{-2}\Lambda(\mathrm dx)\biggr)
\,\mathrm
du \mathop{\sim}\limits_{t \rightarrow0}\frac{t^2}{2-\gamma}
\int_{[1-\exp(-t),1]}x^{-2}\Lambda(\mathrm dx)
\]
and then, applying Karamata's Tauberian theorem (cf. \cite{BGT},
Theorem 1.7.1$'$), that
\vspace*{-3pt}
\begin{eqnarray*}
&&\int_{(0,\infty)} \exp(-nu)u \biggl(
\int_{[1-\exp(-u),1]}x^{-2}\Lambda(\mathrm dx) \biggr)\,\mathrm du
\\[-3pt]
&&\quad
\mathop{\sim}\limits_{n\rightarrow\infty} \frac{\Gamma
(3-\gamma)}{(2-\gamma)n^2}\int_{[1/n,1]}x^{-2}\Lambda(\mathrm dx).
\end{eqnarray*}
Using $\Gamma(3-\gamma)=\Gamma(2-\gamma)(2-\gamma)$, we therefore have,
as $n \rightarrow\infty$
\vspace*{-3pt}
\[
g_n \sim I_n \sim\Gamma(2-\gamma)\int_{[1/n,1]}x^{-2}\Lambda
(\mathrm dx).
\]
\upqed
%\vspace*{-6pt}
\end{pf}
%(here we have used that $x\rightarrow x\exp(-x)$ is bounded on $(0,

\vspace*{-6pt}
\begin{lemma}
\label{lemlamb2} For all measures $\Lambda$ such that
$\int_{[0,1]}x^{-1}\Lambda(\mathrm dx)<\infty$, and all $\lambda
\geq 0$
\vspace*{-3pt}
\[
\sum_{k=1}^{n-1}g_{n,k} \biggl( 1-\biggl(\frac{k}{n}\biggr)^{\lambda} \biggr) \mathop{\rightarrow}\limits_{n\rightarrow\infty}\int_{[0,1]}
\bigl(1-(1-x)^{\lambda}\bigr)x^{-2}\Lambda(\mathrm dx).
\]%\vadjust{\goodbreak}
\end{lemma}
\eject

\begin{pf}
Note that
\begin{eqnarray*}
&&\sum_{k=1}^{n-1}g_{n,k} \biggl( 1-\biggl(\frac{k}{n}\biggr)^{\lambda}\biggr)
\\&&\quad=
\int_{[0,1]}\Biggl(\sum_{k=0}^{n-2} \pmatrix{n\cr k}x^{n-k}(1-x)^{k} \biggl( 1-\biggl(\frac{k+1}{n}\biggr)^{\lambda}\biggr)\Biggr)x^{-2}\Lambda(\mathrm dx)
\\
&&\quad=
\int_{[0,1]}\biggl( \mathbb E\biggl[1-\biggl(\frac{B^{(n,x)}+1}{n}\biggr)^{\lambda}\biggr]-(1-x)^n\biggl(1-\biggl(\frac{n+1}{n}\biggr)^{\lambda}\biggr)\biggr) x^{-2}\Lambda(\mathrm dx),
\end{eqnarray*}
%
%and
%$$ \sum_{k=0}^{n-2} \binom{n}{k} x^{n-k}(1-x)^{k} (
%1-(\frac{k+1}{n})^{\lambda} ) = \mathbb E[
% 1-(\frac{B^{(n,x)}+1}{n})^{\lambda}]-(1-x)^n(1-
%(\frac{n+1}{n})^{\lambda})
%$$
where $B^{(n,x)}$ denotes a binomial random variable with parameters
$n,1-x$. By the strong law of large numbers and dominated convergence
($0 \leq(B^{(n,x)}+1)/n \leq2)$, we have that
\[
\mathbb E\biggl[ 1-\biggl(\frac{B^{(n,x)}+1}{n}\biggr)^{\lambda}\biggr]
\mathop{\rightarrow}\limits_{n \rightarrow\infty}1-(1-x)^{\lambda}\qquad
\forall x
\in[0,1].
\]
Moreover, $(1-x)^n(1-(\frac{n+1}{n})^{\lambda}) \rightarrow0$,
for every $x \in[0,1].$ Besides, since $1-y^{\lambda}
\leq\break\max(1,\lambda)(1-y)$ for $y \in[0,1]$,
\begin{eqnarray*}
\sum_{k=0}^{n-2} \pmatrix{n\cr k} x^{n-k}(1-x)^{k} \biggl(1-\biggl(\frac{k+1}{n}\biggr)^{\lambda} \biggr)
&\leq&
\max(1,\lambda)\sum_{k=0}^{n-2}\pmatrix{n\cr k} x^{n-k}(1-x)^{k} \biggl(1-\biggl(\frac{k+1}{n}\biggr) \biggr)
\\
&\leq&
\max(1,\lambda)\bigl( 1-(1-x+1/n)+(1-x)^{n}/n \bigr)
\\
& \leq&
\max(1,\lambda) x.
\end{eqnarray*}
Using that $\int_{[0,1]}x^{-1}\Lambda(\mathrm dx)<\infty$, we conclude
by dominated convergence.
\end{pf}

\begin{pf*}{Proof of Theorem \ref{theolambda}} Under the
assumptions of Theorem \ref{theolambda}, by Lemmas \ref{lemlamb1} and
\ref{lemlamb2},
\[
1-\sum_{k=0}^n
p_{n,k}\biggl(\frac{k}{n}\biggr)^{\lambda}=\frac{1}{g_n}\Biggl(\sum
_{k=1}^{n-1} g_{n,k}\biggl(1-\biggl(\frac{k}{n}\biggr)^{\lambda}
\biggr) \Biggr)
\sim\frac{\int_{[0,1]} (1-(1-x)^{\lambda})x^{-2}\Lambda(\mathrm
dx)}{\Gamma(2-\gamma) \int_{[1/n,1]} x^{-2}\Lambda(\mathrm d x)}
\]
as $n \rightarrow\infty$ and for all $\lambda\geq0$. Hence (H) holds by Proposition~\ref{sec:preliminaries} and
Theorem \ref{theolambda} is proved.
\end{pf*}

%%%%%%%%%%%%%%%%%%%%%%%%%%%%%%%%
%%%%%%%%%%%%%%%%%%%%%%%%%%%%%%%%

\section*{Acknowledgement}
This work is partially supported by the Agence Nationale
de la Recherche, ANR-08-BLAN-0190 and ANR-08-BLAN-0220-01.

\printhistory

\end{document}